\documentclass[11pt]{article}
\usepackage{latexsym}
\usepackage{amssymb}

\def\qed{\hfill$\Box$}

\newcommand{\Xp}{\mbox{\boldmath $X$}}

\newcommand{\Yp}{\mbox{\boldmath $Y$}}

\newcommand{\orig}{\bf 0}

\newcommand{\weakly}{\mbox{$ \;\stackrel{\cal D}{\longrightarrow}\; $}}

\newcommand{\bydef}{\mbox{$ \;\stackrel{\triangle}{=}\; $}}
\newcommand{\home}{\mbox{$\!\mathtt{\sim}$}}

\newcommand{\approxa}{\mbox{$ \;\stackrel{(a)}{\approx}\; $}}
\newcommand{\approxb}{\mbox{$ \;\stackrel{(b)}{\approx}\; $}}
\newcommand{\approxc}{\mbox{$ \;\stackrel{(c)}{\approx}\; $}}
\newcommand{\approxd}{\mbox{$ \;\stackrel{(d)}{\approx}\; $}}

\newcommand{\RL}{{\mathbb R}}

\newcommand{\IN}{{\mathbb Z}}

\newcommand{\IND}{{\mathbb I}}
\newcommand{\BBP}{{\mathbb P}}
\newcommand{\BBQ}{{\mathbb Q}}
\newcommand{\BBM}{{\mathbb M}}

\newcommand{\PR}{\mbox{\rm Pr}} 
\newcommand{\VAR}{\mbox{\rm Var}} 
\newcommand{\COV}{\mbox{\rm Cov}} 

\newcommand{\essinf}{\mathop{\rm ess\, inf}}
\newcommand{\esssup}{\mathop{\rm ess\, sup}}
\newcommand{\Ahat}{\mbox{$\hat{A}$}}
\newcommand{\Ahatn}{\mbox{$\hat{A}^n$}}

\newcommand{\Ahatnd}{\mbox{$\hat{A}^{n^d}$}}

\newcommand{\Phatn}{\mbox{$\hat{P}_n$}}

\newcommand{\calH}{\mbox{${\cal H}$}}
\newcommand{\calLn}{\mbox{${\cal L}_n$}}

\newcommand{\calR}{\mbox{${\cal R}$}}

\newcommand{\Dmin}{\mbox{$D_{\rm min}$}}
\newcommand{\Dinf}{\mbox{$D_{\rm min}^{(\infty)}$}}

\newcommand{\Dmax}{\mbox{$D_{\rm max}$}}
\newcommand{\Dbar}{\mbox{$\overline{D}$}}

\newcommand{\Dav}{\mbox{$D_{\rm av}$}}
\newcommand{\Dminn}{\mbox{$D_{\rm min}^{(n)}$}}

\newcommand{\LA}{\mbox{$\Lambda$}}

\newcommand{\las}{\mbox{\scriptsize$\lambda$}}
\newcommand{\iid}{\mbox{i.i.d.}\!}
\newcommand{\psipm}{\psi^{\pm}}
\newcommand{\limnd}{\lim_{
	            \mbox{\scriptsize
			 $\begin{array}{c}
				n\to\infty\\
				D\downarrow 0
			  \end{array}$
	     		 }
		 	 }
		   }
\newcommand{\limsupnd}{\limsup_{
                    \mbox{\scriptsize
                         $\begin{array}{c}
                                n\to\infty\\
                                D\downarrow 0
                          \end{array}$
                         }
                         }
                   }
\newcommand{\liminfnd}{\liminf_{
                    \mbox{\scriptsize
                         $\begin{array}{c}
                                n\to\infty\\
                                D\downarrow 0
                          \end{array}$
                         }
                         }
                   }

\newcommand{\la}{\lambda}

\def\be{\begin{eqnarray}}
\def\ee{\end{eqnarray}}
\def\ben{\begin{eqnarray*}}
\def\een{\end{eqnarray*}}

\input{epsf}

\title{Source Coding, Large Deviations,\\
and Approximate Pattern Matching}

\author{A. Dembo \and I. Kontoyiannis}

\date{\today}

\begin{document}
\bibliographystyle{plain}
\maketitle

\thispagestyle{empty}
\setcounter{page}{-2}

\footnotetext[1]{
A. Dembo is with
the Departments of
Mathematics and of
Statistics,
Stanford University,
Stanford, CA 94305.
Email: {\tt amir@stat.stanford.edu}
Web: {\tt www-stat.stanford.edu/\home amir}
}
 
\footnotetext[2]{
I.\ Kontoyiannis is with the Division
of Applied Mathematics, Brown University,
Box F, 182 George St., Providence, RI 02912, USA.
Email: {\tt yiannis@dam.brown.edu}
Web: {\tt www.dam.brown.edu/people/yiannis/}
[Permanent address:
Department of Statistics,
Purdue University,
1399 Mathematical Sciences Building,
W.~Lafayette, IN 47907-1399, USA.]
}
 
\footnotetext[3]{
Amir Dembo was supported in part
by NSF grant \#DMS-0072331.
I. Kontoyiannis was supported in part
by NSF grant \#0073378-CCR.}

\bigskip




\newpage

{\bf Abstract --- }
We present a development of parts of rate-distortion 
theory and pattern-matching algorithms for lossy data 
compression, centered around a lossy version of the 
Asymptotic Equipartition Property (AEP). This treatment 
closely parallels the corresponding development in 
lossless compression, a point of view that was advanced 
in an important paper of Wyner and Ziv in 1989. 
In the lossless case we review how the AEP underlies 
the analysis of the Lempel-Ziv algorithm by viewing it 
as a random code and reducing it to the idealized 
Shannon code. This also provides 
information about the redundancy of the Lempel-Ziv 
algorithm and about the asymptotic behavior of 
several relevant quantities. 

In the lossy case we 
give various versions of the statement of the 
generalized AEP and we outline the general 
methodology of its proof via large deviations. 
Its relationship with Barron and Orey's
generalized AEP is 
also discussed.
The lossy AEP is applied to: (i)~prove strengthened 
versions of Shannon's direct source coding theorem 
and universal coding theorems;
(ii)~characterize the performance of
``mismatched'' codebooks in lossy
data compression;
(iii)~analyze the performance of 
pattern-matching algorithms for lossy 
compression (including Lempel-Ziv schemes); 
(iv)~determine the first order asymptotics 
of waiting times (with distortion) between 
stationary processes; (v)~characterize the 
best achievable rate of ``weighted'' 
codebooks as an optimal sphere-covering 
exponent. We then present a refinement to 
the lossy AEP and use it to: (i)~prove second order 
(direct and converse) lossy source coding theorems,
including universal coding theorems; (ii)~characterize 
which sources are quantitatively easier to compress;
(iii)~determine the second order
asymptotics of waiting times between
stationary processes;
(iv)~determine the precise asymptotic
behavior of longest match-lengths 
between stationary processes.
Extensions to random fields are also given.

\medskip

{\bf Index Terms --- } Rate-distortion theory, 
pattern-matching, large deviations, 
data compression.

\bigskip

\tableofcontents

\newpage 
\section{Introduction}

\subsection{Lossless Data Compression}
It is probably only a slight exaggeration to 
say that the central piece of mathematics
in the proof of almost any lossless coding
theorem is provided by the
Asymptotic Equipartition Property, or AEP.
Suppose we want to (losslessly) 
compress a message 
$X_1^n=(X_1,X_2,\ldots,X_n)$
generated by a stationary 
memoryless source $\Xp=\{X_n\;;\;n\geq 1\}$
where each $X_i$ takes values in the 
finite alphabet $A$
(much more general situations will be considered
later). For this source, the AEP states that
as $n\to\infty$
\be
-\frac{1}{n}\log_2 P^n(X_1^n)\to H
\;\;\;\;\mbox{in probability}
\label{eq:shannonAEP}
\ee
where $P$ is the common distribution of the 
independent and identically distributed ($\iid$) 
random variables $X_i$, $P^n$ denotes the 
(product) joint distribution of $X_1^n$, and 
$H=E[-\log_2 P(X_1)]$ is
the entropy rate of the source --
see Shannon's original paper
\cite[Theorem~3]{shannon:48} or 
Cover and Thomas' text \cite[Chapter~4]{cover:book}.
[Here and throughout the paper, 
$\log_2$ denotes the logarithm taken
to base 2, and $\log$ denotes
the natural logarithm.]
From (\ref{eq:shannonAEP}) we can 
immediately extract some useful 
information: It implies that 
when $n$ is large the message 
$X_1^n$ will most likely have 
probability at least as high 
as $2^{-n(H+\epsilon)}$:
\be
P^n(X_1^n)\geq 2^{-n(H+\epsilon)}
	\;\;\;\;\mbox{with high probability.}
\label{eq:tocompare}
\ee
But there cannot be many high-probability messages.
In fact, there can be at most $2^{n(H+\epsilon)}$
messages with $P^n(X_1^n)\geq 2^{-n(H+\epsilon)}$,
so we need approximately $2^{nH}$ representative
messages from the source $\Xp$ in order to cover
our bets (with high probability). 
If we let ${\cal T}_n$ be the 
set of high-probability
strings $x_1^n\in A^n$ having 
$P^n(x_1^n)\geq 2^{-n(H+\epsilon)}$,
then with high probability we can
correctly represent the source output
$X_1^n$ by an element of ${\cal T}_n$.
Since there are no more than 
$2^{n(H+\epsilon)}$ of them, 
we need no more than $nH$ bits
to correctly encode $X_1^n$.

\paragraph{Shannon's Random Code.}
Another way to extract information
from (\ref{eq:shannonAEP}) is as follows. 
The fact that for large $n$ we typically have
$P^n(X_1^n)\approx 2^{-nH}$ also means that
if we independently generate another random
string, say $Y_1^n$, from the same distribution
as the source, the probability that $X_1^n$
is the same as $Y_1^n$ is about $2^{-nH}$.
Suppose that instead of using the strings in
${\cal T}_n$ above as our representatives
for the source, we decided to independently
generate a collection of random strings 
$Y_1^n$ from the distribution $P^n$; how
many would we need? Given a source string
$X_1^n$, the probability that any one
of the $Y_1^n$ matches it is 
$\approx 2^{-nH}$, so in order to
have high probability of success
in representing $X_1^n$ without error
we should choose approximately
$2^{n(H+\epsilon)}$ random strings $Y_1^n$.
Therefore, whether we choose the set of
representatives systematically or
randomly, we always need about
$2^{nH}$ strings in order to be able
to encode $X_1^n$ losslessly with high
probability. Note that the randomly
generated set ${\cal T}_n$ is nothing
but Shannon's random codebook
\cite{shannon:59} specialized to the
case of lossless compression.

\paragraph{Idealized Lempel-Ziv Coding.}
In 1989, in a very influential paper 
\cite{wyner-ziv:1}, Wyner and Ziv took 
the above argument several steps further. 
Aiming to ``obtain 
insight into the workings of [...] the 
Lempel-Ziv data compression algorithm,'' 
they considered the following coding 
scenario: Suppose that an encoder and a
decoder both have available to them
a long database, say an infinitely 
long string $Y_1^\infty=(Y_1,Y_2,\ldots)$
that is independently generated from 
the same distribution as the source. 
Given a source string $X_1^n$ to 
be transmitted, the
encoder looks for the
first appearance of $X_1^n$ in the
database (assuming, for now,
that it does appear somewhere).
Let $W$ denote the position of
this first appearance, that is,
let $W$ be the smallest integer
for which 
$Y_W^{W+n-1}=(Y_W,Y_{W+1},\ldots,Y_{W+n-1})$
is equal to $X_1^n$.
Then all the encoder has to do is 
it to tell the decoder the value of 
$W$; the decoder can read off the string 
$Y_W^{W+n-1}$ and recover $X_1^n$ 
perfectly. This description can be
given using
(cf. \cite{elias}\cite{wyner-ziv:2})
no more than 
\be
\ell(X_1^n)=\log_2 W + O(\log_2\log_2 W)
	\;\;\;\;\mbox{bits.}
\label{eq:elias1}
\ee

How good is this scheme? 
First note that, for any given source
string $X_1^n$, the random variable
$W$ records the first ``success'' in a 
sequence of trials (``Is $Y_1^n=X_1^n$?,''
``Is $Y_2^{n+1}=X_1^n$?,''
and so on),
each of which has probability
of success $p=P^n(X_1^n)$. Although
these trials are not independent,
for large $n$ they are almost independent
(in a sense that 
will be made precise below), so the 
distribution of $W$ is close to 
a geometric with parameter 
$p=P^n(X_1^n)$.
For long strings $X_1^n$ (i.e., for
large $n$) $p$ is small,
and $W$ is typically close to its 
expected value, which is approximately 
equal to the mean of a geometric
random variable
with parameter $p$, namely $1/p$. But the 
AEP tells us that, when $n$ is large,
$p=P^n(X_1^n)\approx 2^{-nH}$, so we 
expect $W$ to  be typically around $2^{nH}$.
Hence, from (\ref{eq:elias1}) the 
description length $\ell(X_1^n)$
of $X_1^n$ will be, to first order,
$$\ell(X_1^n)
\approx -\log_2 P^n(X_1^n) 
\approx nH 
	\;\;\;\mbox{bits, with high probability.}$$
This shows that above scheme is asymptotically
optimal, in that its limiting compression ratio 
is equal to the entropy.

\paragraph{Practical Lempel-Ziv Coding.}
The Lempel-Ziv algorithm 
\cite{ziv-lempel:1}\cite{ziv-lempel:2} and
its many variants (see, e.g.,  
\cite[Ch.~8]{bell:cleary:witten}) are some
of the most successful data compression
algorithms used in practice. Roughly speaking, 
the main idea behind these algorithms is to 
use the message's own past as a database 
for future encoding. Instead of looking
for the first match in an infinitely long
database, in practice the encoder looks
for the longest match in a database of
fixed length. The analysis in 
\cite{wyner-ziv:1} of the idealized 
scheme described above was the first 
step in providing a probabilistic 
justification for the optimality
of the actual practical algorithms. 
Subsequently, in \cite{wyner-ziv:3} 
and \cite{wyner-ziv:2} Wyner and Ziv
established the asymptotic optimality
of the Sliding-Window (SWLZ) and
the Fixed-Database (FDLZ) versions
of the algorithm.

\subsection{Lossy Data Compression}

A similar development to the one outlined
above can be given in the case of lossy
data compression, this time centered around
a lossy analog of the AEP \cite{my:thesis}. 
To motivate this
discussion we look at Shannon's original 
random coding proof of the (direct) lossy 
source coding theorem \cite{shannon:59}. 

\paragraph{Shannon's Random Code.}
Suppose we want to describe the output 
$X_1^n$ of a memoryless source,
with distortion $D$ or less with respect 
to a family of single-letter distortion 
measures $\{\rho_n\}$. 
Let $Q_n^*$ be the optimum reproduction 
distribution on $\Ahatn$, where 
$\Ahat$ is the reproduction alphabet. 
Shannon's random coding 
argument says that we should 
construct a codebook ${\cal T}_n$ 
of $2^{n(R(D)+\epsilon)}$ codewords 
$Y_1^n$ generated $\iid$ from $Q_n^*$,
where $R(D)$ is the rate-distortion
function of the source (in bits).
The proof that $2^{n(R(D)+\epsilon)}$ codewords 
indeed suffice is based on the following 
result, Lemma~1 in \cite{shannon:59}.

\medskip

{\em Shannon's ``Lemma~1'':}
For $x_1^n\in A^n$ let $B(x_1^n,D)$ denote 
the distortion-ball of radius $D$ around
$x_1^n$, i.e., the collection of all
reproduction strings $y_1^n\in\Ahatn$ with
$\rho_n(x_1^n,y_1^n)\leq D$.
When $n$ is large:\footnote{The
notation in Shannon's statement is
slightly different, and he considers
the more general case of ergodic sources. 
For the sake of clarity we restrict
attention here to the $\iid$ case.}
\be
Q_n^*(B(X_1^n,D))\geq 2^{-n(R(D)+\epsilon)}
	\;\;\;\;\mbox{with high probability.}
\label{eq:lemma1}
\ee

\medskip

In the proof of the coding theorem
this 
lemma plays the same role that
the AEP played in the lossless case;
notice the similarity between 
(\ref{eq:lemma1}) and its analog 
(\ref{eq:tocompare}) in the lossless case. 

Let's fix a source string $X_1^n$ to 
be encoded. The probability that $X_1^n$ 
matches any one of the codewords 
$Y_1^n$ in ${\cal T}_n$ is
$$\Pr\{\rho_n(X_1^n,Y_1^n)\leq D\,|\,X_1^n\}
=
\Pr\{Y_1^n\in B(X_1^n,D)\,|\,X_1^n\}
=
Q_n^*(B(X_1^n,D))$$
and by the lemma this probability is 
at least $2^{-n(R(D)+\epsilon)}$.
Therefore, with $2^{n(R(D)+\epsilon)}$
independent codewords to choose from,
we have a good chance for finding
a match with distortion $D$ or less.

\paragraph{Generalized AEP and Applications.}
A stronger and more general version
of Lemma~1 will be our starting point
in this paper. In the following section 
we will prove a {\em generalized AEP}:
For any product measure 
$Q^n$ on $\Ahatn$
\be
-\frac{1}{n}\log Q^n(B(X_1^n,D)) \to R_1(P,Q,D)
	\;\;\;\;\mbox{w.p.1}
\label{eq:firstDAEP}
\ee
where 
$R_1(P,Q,D)$ is a (non-random) function
of the distributions $P$ and $Q$ and
of the distortion level $D$. 
[We will later prove several
variants of (\ref{eq:firstDAEP})
under much weaker assumptions.]

Like the AEP in the lossless case,
the generalized AEP and its refinements
find numerous applications in data 
compression, universal data compression, 
and in general pattern-matching questions.
Many of these applications were inspired
by the treatment in Wyner and Ziv's 1989
paper \cite{wyner-ziv:1}.  A (very
incomplete) sample of subsequent
work in the Wyner-Ziv spirit
includes the papers
	\cite{steinberg-gutman}\cite{luczak-szpankowski}\cite{
	yang-kieffer:1}\cite{kontoyiannis-lossy1-1}
	on lossy data compression, 
	and
	\cite{luczak-szpankowski}\cite{
	dembo-kontoyiannis}\cite{yang-zhang:99c}
	on pattern-matching.

Aaron Wyner himself remained active in this
field for the following ten years, and his
last paper \cite{wyner-ziv-wyner}, co-written
with J.~Ziv and A.J.~Wyner,
was a review paper on this subject.
In the present paper we review the corresponding
developments in the lossy case, and in the process
we add new results (and some new proofs of 
recent results) in an attempt to present a more
complete picture.

\subsection{Central Themes, Paper Outline}

In Section~2 we give an extensive discussion of
the generalized AEP. By now there are numerous
different proofs under different assumptions, 
and we offer a streamlined approach to the most 
general versions using techniques from large 
deviation theory (cf.
\cite{yang-kieffer:1}\cite{dembo-kontoyiannis}%
\cite{chi-it:01}\cite{chi-AP:01}
and Bucklew's earlier work 
\cite{bucklew:87}\cite{bucklew:88}). We also
discuss the relationship 
of the generalized AEP
with the classical extensions
of the AEP (due to Barron \cite{barron:1}
and Orey \cite{orey:85})
to processes with densities.
We establish a formal connection
between these two by looking at the limit
of the distortion level $D\downarrow 0$.

In Section~3 we develop applications 
of the generalized AEP to a number
of related problems. 
We show how the generalized AEP
can be used to determine the
asymptotic behavior of Shannon's
random coding scheme, and we
discuss the role of mismatch
in
lossy data compression.
We also determine the first order 
asymptotic behavior of 
waiting times and longest 
match-lengths between stationary 
processes.  The main ideas used 
here are strong approximation 
\cite{kontoyiannis-jtp} and 
duality \cite{wyner-ziv:1}.
We present strengthened versions
of Shannon's direct lossy source coding
theorem (and of a corresponding universal
coding theorem), showing that {\em almost all}
random codebooks achieve essentially 
the same compression performance. 
A lossy version of the Lempel-Ziv 
algorithm is recalled, which 
achieves optimal compression 
performance (asymptotically)
as well as polynomial 
complexity at the encoder.
We also discuss how the classical 
source coding problem
can be generalized to a question about
weighted sphere-covering. The answer 
to this question gives, as
corollaries, Shannon's coding theorems,
Stein's lemma in hypothesis testing, 
and some converse concentration inequalities.

Section~4 is devoted to second order
refinements of the AEP and the generalized
AEP. It is shown, for example, that 
under certain conditions
$-\log P^n(X_1^n)$ and $-\log Q^n(B(X_1^n,D))$ 
are asymptotically Gaussian.
These refinements are used in Section~5
to provide corresponding second order 
results (such as central limit theorems) 
for the applications considered in Section~3. 
We prove second order asymptotic results 
for waiting times
and longest match-lengths. 
Precise redundancy rates are 
given for Shannon's random code,
and converse coding theorems show 
that the random code achieves the 
optimal pointwise redundancy, 
up to terms of order $(\log n)$. 
For $\iid$ sources the pointwise 
redundancy is typically of order 
$\sigma\sqrt{n}$, where $\sigma$
is the minimal coding variance of 
the source. When $\sigma=0$ these
fluctuations disappear, and the
best pointwise redundancy is of
order $(\log n)$. The question of 
exactly when $\sigma$ can be equal 
to zero is briefly discussed.

Finally, Sections~6 and 7 
contain generalizations of some
of the above results to 
random fields. All the results
stated there
are new, although
most of them are straightforward 
generalizations of corresponding
one-dimensional results.

\section{The Generalized AEP}

\subsection{Notation and Definitions}
We begin by introducing some basic definitions
and notation that will remain in effect for 
the rest of the paper. We will consider a
stationary ergodic 
process 
$\Xp=\{X_n\;;\;n\in\IN\}$
taking values in a general alphabet 
$A$.\footnote{To avoid 
uninteresting technicalities,
we will assume throughout that 
$A$ is a complete, separable metric space,
equipped with its associated Borel 
$\sigma$-field ${\cal A}$.
Similarly we take $(\Ahat,\hat{\cal A})$
to be the Borel measurable space
corresponding to a complete,
separable metric space $\Ahat$.} When
talking about data compression, $\Xp$
will be our source and $A$ will be 
called the source alphabet. We write
$X_i^j$ for the vector of random
variables $X_i^j=(X_i,X_{i+1},\ldots,X_j)$,
and similarly 
$x_i^j=(x_i,x_{i+1},\ldots,x_j)\in 
A^{j-i+1}$
for a realization of these random variables,
$-\infty\leq i\leq j\leq \infty$.
We let $P_n$ denote the marginal 
distribution of $X_1^n$ on $A^n$
($n\geq 1$), and write $\BBP$
for the distribution of the whole
process.
Similarly, we take
$\Yp=\{Y_n\;;\;n\in\IN\}$
to be a stationary ergodic
process taking values in the
(possibly different) alphabet
$\Ahat$.${}^2$ 
In the context of
data compression,
$\Ahat$ is the reproduction
alphabet and $\Yp$
has the ``codebook'' distribution.
We write $Q_n$ for the marginal
distribution of $Y_1^n$ on $\Ahatn$,
$n\geq 1$, and $\BBQ$ for the
distribution of the whole process $\Yp$.
We will always assume that the process
$\Yp$ is independent of $\Xp$.

Let $\rho:A\times\Ahat\to[0,\infty)$
be an arbitrary nonnegative (measurable)
function, and define a sequence of
single-letter distortion measures 
$\rho_n:A^n\times\Ahatn\to[0,\infty)$ by
\ben
\rho_n(x_1^n,y_1^n)\bydef\frac{1}{n}\sum_{i=1}^n\rho(x_i,y_i)
\;\;\;\;x_1^n\in A^n,\;y_1^n\in\Ahatn.
\een
Given $D\geq 0$ and $x_1^n\in A^n$, 
we write 
$B(x_1^n,D)$ for
the distortion-ball of radius $D$ around $x_1^n$:
$$B(x_1^n,D)=\{y_1^n\in\Ahatn\;:\;\rho_n(x_1^n,y_1^n)\leq D\}.$$

Throughout the paper, $\log$ denotes
the natural logarithm and $\log_2$
the logarithm to base 2. Unless otherwise
mentioned, all familiar information-theoretic
quantities (such as the entropy,
mutual information, and so on)
are assume to be defined in terms
of natural logarithms (and are 
therefore given in nats).

\subsection{Generalized AEP When $\Yp$ is I.I.D.}

In the case when $A$ is finite,
the classical AEP, also known as the 
Shannon-McMillan-Breiman theorem
(see \cite[Chapter~15]{cover:book} 
or the original papers
\cite{shannon:48}\cite{mcmillan}\cite{breiman:57}\cite{breiman:60}),
states that as $n\to\infty$
\be
-\frac{1}{n}\log P_n(X_1^n) \to H(\BBP)
\;\;\;\;\mbox{w.p.1}
\label{eq:discreteAEP}
\ee
where 
$$H(\BBP)\bydef\lim_{n\to\infty}\frac{1}{n}H(X_1^n)$$
is the entropy rate of the process $\Xp$
(in nats, since we are taking logarithms to
base $e$).
As we saw in the
Introduction,
in lossy data
compression the role of the AEP is taken 
up by the result of Shannon's ``Lemma 1'' 
and, more generally, by statements of the form
$$-\frac{1}{n}Q_n(B(X_1^n,D))\to R(\BBP,\BBQ,D)
\;\;\;\;\mbox{w.p.1}$$
for some 
non-random
``rate-function'' $R(\BBP,\BBQ,D)$.

First we consider the simplest case where $\Yp$ 
is assumed to be an $\iid$ process. We write 
$Q=Q_1$ for its first order marginal,
so that $Q_n=Q^n$, for $n\geq 1.$ Similarly
we write $P=P_1$ for the first order marginal 
of $\Xp$. 
Let 
\be
\Dmin & \bydef &  E_P[\essinf_{Y\sim Q} \;\rho(X,Y)]
	\label{eq:Dmin} \\
\Dav & \bydef &  E_{P\times Q}[\rho(X,Y)].
	\label{eq:Dav}
\ee
[Recall that the essential infimum of
a function $g(Y)$ of the random variable
$Y$ with distribution $Q$ is defined as
$\essinf_{Y\sim Q} g(Y) = 
\sup\{t\in\RL\;:\;Q\{g(Y)>t\}=1\}.$]

Clearly $0\leq\Dmin\leq\Dav$.
To avoid the trivial case when $\rho(x,y)$
is essentially constant for
($\BBP$-almost) all
$x\in A$, we assume that with positive
$\BBP$-probability $\rho(x,y)$ is not
essentially constant in $y$, that is:
\be
\Dmin < \Dav.
\label{eq:nonconst}
\ee
Note also that
for $D$ greater than $\Dav$,
the probability 
$Q^n(B(X_1^n,D))\to 1$
as $n\to\infty$
(this is easy to see 
by the ergodic theorem),
so we restrict our attention
to distortion levels $D<\Dav$.

\medskip 

{\em Theorem~1. Generalized AEP when $\Yp$ is $\iid$:}
Let $\Xp$ be a stationary ergodic process
and $\Yp$ be $\iid$ with marginal distribution 
$Q$ on $\Ahat$.
Assume that $\Dav=E_{P\times Q}[\rho(X,Y)]$ is
finite. Then for any $D\in(\Dmin,\Dav)$
\ben
-\frac{1}{n}\log Q^n(B(X_1^n,D)) \to R_1(P,Q,D)
        \;\;\;\;\mbox{w.p.1}.
\een
The rate-function $R_1(P,Q,D)$ 
is defined as
\ben
R_1(P,Q,D) = \inf_W H(W\|P\times Q)
\een
where $H(W\|V)$ denotes the relative
entropy between two distributions
$W$ and $V$,
$$H(W\|V) \bydef \left\{ \begin{array}{ll}
   E_W[\log\frac{dW}{dV}]
		   & \mbox{if the density $\frac{dW}{dV}$ exists}, \\
   \infty 	   & \mbox{otherwise}
 \end{array} \right.
$$
and the infimum is taken over all 
joint distributions $W$ on
$A\times\Ahat$ such that 
the first marginal of $W$ is $P$
and $E_W[\rho(X,Y)]\leq D.$

\medskip

{\em Example~1: The rate-function $R_1(P,Q,D)$ 
when $Q$ is Gaussian:} 
Although in general the rate-function
$R_1(P,Q,D)$ cannot be evaluated explicitly,
here we show that it is possible to obtain 
an exact expression for $R_1(P,Q,D)$ in the 
special case when $\rho(x,y)=(x-y)^2$,
$\Xp$ is a real-valued, process, 
and $Q$ is a Gaussian measure 
on $\RL.$ Specifically, assume
that $\Xp$ is a zero-mean,
stationary ergodic process
with finite variance 
$\sigma^2=\VAR(X_1)<\infty$,
and take $Q$ to be 
a zero-mean Gaussian measure 
with variance $\tau^2$, i.e.,
$Q\sim N(0,\tau^2)$.
Under these assumptions, it is easy to see
that $\Dmin=0$ and $\Dav=\sigma^2+\tau^2$.
Moreover, with the help of Proposition~2 below,
$R_1(P,Q,D)$ can be explicitly
evaluated as:
$$R_1(P,Q,D) = \left\{ \begin{array}{ll}
	\infty\,,	& \;\;\;\;D=0\\
	\frac{1}{2}\log\left(\frac{v}{D}\right)
	  -\frac{(v-D)(v-\sigma^2)}
		{2v\tau^2}\,,
	 		& \;\;\;\;0<D<\sigma^2+\tau^2\\
	0\,, 		& \;\;\;\;D\geq \sigma^2+\tau^2
 \end{array} \right.
$$
where 
$$v\bydef\frac{1}{2}\left[\tau^2+\sqrt{\tau^4+4D\sigma^2}\right].$$
We will come back to this example when considering
mismatched rate-distortion codebooks in Section~3.2.

\medskip
 
{\em Remark 1:}
In more familiar information-theoretic
terms, the rate-function $R_1(P,Q,D)$
can equivalently be defined as 
(cf. \cite{yang-kieffer:1})
\ben
R_1(P,Q,D) = \inf_{(X,Y)}\,[I(X;Y)
	+H(Q_Y\|Q)]
\een
where $I(X;Y)$ denotes the mutual
information (in nats) between the 
random variables $X$ and $Y$,
and the infimum is over
all jointly distributed random variables
$(X,Y)$ with values in $A\times\Ahat$
such that $X$ has distribution $P$,
$E[\rho(X,Y)]\leq D$, and $Q_Y$ denotes
the distribution of $Y$.

\medskip

{\em Remark 2:}
The assumption that $\Yp$ is
$\iid$ is clearly restrictive
and it will be relaxed below.
On the other hand the assumptions
on the distortion measure 
$\rho$ seem to be minimal;
we simply assume that $\rho$
has finite expectation (in
the more general results below
$\rho$ is assumed to be bounded).
In this form, the result of 
Theorem~1 is new.

\medskip

{\em Discussion of Proof:}
Let's fix a realization $x_1^\infty$
of $\Xp$. The probability
$Q^n(B(X_1^n,D))$ can be written as
$$\PR\left\{Y_1^n\in B(X_1^n,D) \,|\,X_1^n=x_1^n\right\} \;=\;
\PR\left\{\frac{1}{n}\sum_{i=1}^n\rho(x_i,Y_i)\leq D \right\}.$$
Since the distortion level $D$ is taken smaller than the
average value $\Dav$, this is large deviations probability
for the partial sums $(1/n)\sum_{i=1}^n Z_i$ of
the independent (but not identically distributed)
random variables $Z_i=\rho(x_i,Y_i)$. The proof
is essentially an application of the 
G\"artner-Ellis theorem of large deviations
to the random variables $\{Z_i\}$.

\medskip

{\em Proof Outline:}
Choose and fix a realization $x_1^\infty$ of
$\Xp$ and define the random variables 
$Z_i=\rho(x_i,Y_i)$. Let 
$$S_n=\frac{1}{n}\sum_{i=1}^nZ_i$$
and define the log-moment generating
functions of the normalized partial 
sums $S_n$ by
$$\LA_n(\la) \bydef \log
	E_{Q^n}\left(e^{\lambda S_n}\right),
	\;\;\;\;\lambda\leq 0.$$
Then for any $\la\leq 0$, by the ergodic theorem we have
that
\be
\frac{1}{n}\LA_n(n\la)
	= \frac{1}{n}\sum_{i=1}^n\log
	E_Q\left(e^{\lambda\rho(x_i,Y_i)}\right)
	\to 
	\LA(\la)\bydef E_P\left[\log
	E_Q\left(e^{\lambda\rho(X,Y)}\right)
	\right]
\label{eq:GEcheck}
\ee
for $\BBP$-almost any realization $x_1^\infty$.
Now we would like to apply the 
G\"artner-Ellis theorem, but first
we need to check some simple properties
of the function $\LA(\la)$. 
Note that 
$\LA(\la)\leq 0$ 
and
also (by Jensen's inequality)
$\LA(\la)\geq \la\Dav>-\infty$,
for all 
$\la\leq 0$.
Moreover, $\LA(\la)$ is twice
differentiable
in $\la$ with
$$\LA'(\la) = E_{P\times Q}\left(\rho(X,Y)
    \frac
	{e^{\lambda\rho(X,Y)}}
	{E_Q[e^{\lambda\rho(X,Y')}]}
			   \right)
$$
and 
$$\LA''(\la) = 
E_P\left[
    E_Q
    \left\{\rho^2(X,Y)
    \frac
        {e^{\lambda\rho(X,Y)}}
        {E_Q[e^{\lambda\rho(X,Y')}]}
    \right\}
	\;-\;
    \left(E_Q
        \left\{
	\rho(X,Y)
    	\frac
        {e^{\lambda\rho(X,Y)}}
        {E_Q[e^{\lambda\rho(X,Y')}]}
        \right\}
    \right)^2
\right]$$
(this differentiability is easily verified by
an application of the dominated convergence
theorem). By the Cauchy-Schwarz
inequality $\LA''(\la)\geq 0$ for all $\la<0$,
and in fact $\LA''(\la)$ is strictly positive
due to assumption (\ref{eq:nonconst}).
Also it is not hard to verify that 
$$\lim_{\lambda\uparrow 0}\LA'(\la)=\Dav$$
and
\be
\lim_{\lambda\downarrow -\infty}\LA'(\la)=\Dmin.
\label{eq:la-lim}
\ee
Since $D\in(\Dmin,\Dav)$,
there exists a unique $\la^*<0$ with
$\LA'(\la^*)=D$, and therefore
the Fenchel-Legendre
transform of $\LA(\la)$ evaluated at $D$ is
$$\LA^*(D)\bydef\sup_{\la\leq 0}[\la D-\LA(\la)]
	\;=\;\la^*D-\LA(\la^*).$$
Now we can apply the
G\"artner-Ellis theorem
\cite[Theorem~2.3.6]{dembo-zeitouni:book}
to deduce from
(\ref{eq:GEcheck})
that with $\BBP$-probability one
$$-\frac{1}{n}\log Q^n(B(X_1^n,D)) \to \LA^*(D).$$
The proof is complete upon noticing 
that $\LA^*(D)$ is nothing but $R_1(P,Q,D)$.
This is stated and proved in 
the following proposition.
\qed

\medskip

{\em Proposition 2. Characterization of the Rate Function:}
In the notation of the proof of Theorem~1,
$\LA^*(D)=R_1(P,Q,D)$, for $D\in(\Dmin,\Dav)$.

\medskip

{\em Proof Outline:} Under additional 
assumptions on the distortion measure 
$\rho$ this has appeared in various papers
(see, e.g., \cite{dembo-kontoyiannis}\cite{yang-zhang:99}).
For completeness, we offer a proof sketch here.

In the notation of the above proof, consider
the measure $W$ on $A\times\Ahat$ defined by
$$\frac{dW(x,y)}{dP\times Q} = 
\frac{e^{\las^*\rho(x,y)}}
{E_Q[e^{
\las^* 
\rho(x,Y)}]}.$$
Obviously the first marginal of $W$
is $P$ and it is easy to check that 
that $E_W[\rho(X,Y)]=\LA'(\la^*)=D$.
Therefore, by the definitions of
$R_1(P,Q,D)$ and $W$, and by
the choice of $\la^*$:
\be
R_1(P,Q,D)\leq H(W\|P\times Q)
	=\la^*D-\LA(\la^*)
	=\LA^*(D).
\label{eq:propUBD}
\ee
To prove the corresponding lower
bound we first claim that 
for any measurable function
$\phi:\Ahat\to (-\infty,0]$,
and any probability measure 
$Q'$ on $\Ahat$,
\be
H(Q'\|Q)\geq E_{Q'}(\phi(Y)) -\log E_{Q}(e^{\phi(Y)}).
\label{eq:generalSV}
\ee
Let $Q_\phi$ denote the probability measure on $\Ahat$ such that
$dQ_\phi/dQ=e^\phi/E_{Q}(e^{\phi(Y)})$. Clearly, it 
suffices to prove (\ref{eq:generalSV}) in case $dQ'/dQ$ exists,
in which case the difference between the left and right hand sides is 
$$
 E_{Q'}\left\{\log\frac{dQ'}{dQ}\right\} -
 E_{Q'}\left\{\log\left(\frac{e^{\phi}}
	{
	E_{Q}(e^{\phi})
	}
	\right)\right\}
\;=\; 
%
H(Q'\|
Q_\phi)
\;\geq\;0.
$$
Given an arbitrary
candidate $W$ as in the definition of
$R_1(P,Q,D)$ and any $x\in A$, we take
$Q'=W(\cdot|x)$ and $\phi(y)=\la^*\rho(x,y)$
in (\ref{eq:generalSV}) to get that 
$$H(W(\cdot|x)\|Q(\cdot))\geq \la^* E_{W(Y|x)}[\rho(x,Y)]
	- \log E_{Q}(e^{\lambda^*\rho(x,y)}).$$
Substituting $X$ for $x$,
taking expectations of both sides 
with respect to $P$, 
and recalling that $\la^*<0$
and $E_W[\rho(X,Y)]\leq D$, we get:
$$H(W\|Q)\geq \la^*D - \LA(\la^*) = \LA^*(D).$$
Since $W$ was arbitrary it follows that
$R_1(P,Q,D)\geq \LA^*(D)$, and together
with (\ref{eq:propUBD}) this completes 
the proof.
\qed

\subsection{Generalized AEP When $\Yp$ is Not I.I.D.}

Next we present two versions 
of the generalized AEP that hold
when $\Yp$ is a stationary 
dependent process,
under some additional conditions.

Throughout this section we will
assume that the distortion
measure is essentially bounded
\be
\Dmax \bydef \esssup_{(X_1,Y_1)\sim P_1\times Q_1} 
		\rho(X_1,Y_1)<\infty.
\label{eq:Dmax}
\ee
We let $\Dav$ be defined as earlier,
$\Dav = E_{P_1\times Q_1}[\rho(X_1,Y_1)]$,
and for $n\geq 1$ we let
\ben
\Dminn 
\bydef
	E_{P_n}
	\left[
		\essinf_{Y_1^n\sim Q_n} \;\rho_n(X_1^n,Y_1^n)
	\right].
\een
It is easy to see that $n \Dminn$ is a finite, superadditive sequence, 
and therefore we can also define
$$\Dmin = \lim_{n\to\infty} \Dminn = \sup_{n\geq 1} \Dminn.$$
As before, we will assume that 
the distortion measure $\rho$ is 
not essentially constant,
that is, $\Dmin<\Dav.$

\medskip

We first state a version 
of the generalized AEP that 
was recently proved by Chi 
\cite{chi-it:01}, for 
processes $\Yp$ satisfying 
a rather strong
mixing condition: We say that 
the stationary process $\Yp$ 
is {\em $\psipm$-mixing}, if 
for all $d$ large enough there is a
finite constant $c_d$ such that
$$
c_d^{-1} \BBQ(A)\BBQ(B) < \BBQ(A\cap B) < c_d \BBQ(A)\BBQ(B)
$$
for all events $A\in\sigma(Y_{-\infty}^0)$
and $B\in\sigma(Y_d^{\infty})$,
where $\sigma(Y_i^j)$ denotes
the $\sigma$-field generated by $Y_i^j$. 
Recall the usual definition 
according to which $\Yp$ is called 
{\em $\psi$-mixing} if in fact 
the constants $c_d \to 1$ as 
$d \to \infty$; see \cite{bradley} 
for more details.
Clearly $\psipm$-mixing is weaker 
than {\em $\psi$-mixing}.

\medskip


{\em Theorem~3. Generalized AEP when $\Yp$ is $\psipm$-mixing
\cite{chi-it:01}:}
Let $\Xp$ and $\Yp$ be stationary ergodic
processes. Assume that $\Yp$ is $\psipm$-mixing,
and that the distortion measure $\rho$ is bounded.
Then for all $D\in(\Dmin,\Dav)$
\be
-\frac{1}{n}\log Q_n(B(X_1^n,D)) \to R(\BBP,\BBQ,D)
        \;\;\;\;\mbox{w.p.1}
\label{eq:thm4}
\ee
where $R(\BBP,\BBQ,D)$ is the rate-function defined by
\be
\label{eq:thm4b}
R(\BBP,\BBQ,D) = \lim_{n\to\infty} R_n(P_n,Q_n,D) 
\ee
where, for $n \geq 1$, 
\ben
R_n(P_n,Q_n,D) \bydef \inf_{V_n} n^{-1} H(V_n\|P_n\times Q_n)
\een
and the infimum is taken over all joint 
distributions $V_n$ on $A^n\times\Ahatn$ 
such that the $A^n$-marginal of $V_n$ 
is $P_n$ and $E_{V_n}[\rho_n(X_1^n,Y_1^n)]\leq D$.

\medskip

As we discussed in the previous section,
the proof of most versions of the generalized
AEP consistst of two steps: First a 
``conditional large deviations'' result
is proved for the random variables 
$\{\rho_n(x_1^n,Y_1^n)\;;\;n\geq 1\}$,
where $x_1^\infty$ is a fixed realization of the
process $\Xp$. Second, the rate-function
$R(\BBP,\BBQ,D)$ is characterized as the 
limit of a sequence of minimizations in 
terms of relative entropy.

In a subseqeunt paper, Chi \cite{chi-AP:01}
showed that the first of these steps 
(the large deviations part) remains
valid under a condition 
weaker than $\psipm$-mixing, 
condition~$(S)$ of \cite{bryc-dembo:96}.
In the following theorem we give a general
version of the second step; we prove
that the generalized AEP (\ref{eq:thm4}) 
and the formula (\ref{eq:thm4b}) for the 
rate-function remain valid as long as 
the random variables 
$\{\rho_n(x_1^n,Y_1^n)\;;\;n\geq 1\}$
satisfy a large deviations principle (LDP)
with some {\it deterministic}, convex 
rate-function (see \cite{dembo-zeitouni:book}
for the precise meaning of this statement).

\medskip

{\em Theorem~4.} 
Let $\Xp$ and $\Yp$ be stationary processes.
Assume that $\rho$ is bounded, and that with
$\BBP$-probability one, conditional on 
$X_1^\infty=x_1^\infty$, the random variables
$\{\rho_n(x_1^n,Y_1^n)\;;\;n\geq 1\}$ satisfy a
large deviations principle with some
deterministic, convex rate-function.
Then, both (\ref{eq:thm4}) and (\ref{eq:thm4b})
hold for any $D\in(\Dmin,\Dav)$, except 
possibly at the point $D=\Dinf$, where
\be
\Dinf \bydef \inf \{ D \geq 0 : 
\sup_{n \geq 1} R_n(P_n,Q_n,D) < \infty \}.
\label{eq:dinf}
\ee

\medskip

Since Theorem~4 has an exact 
analog in the case of random fields, 
we postpone its proof until the
proof of the corresponding result
(Theorem~27) in Section~6.

\medskip

{\em Remark 3:} 
Suppose that the joint process $(\Xp,\Yp)$ is
stationary, and that it satisfies a
``process-level large deviations principle''
(see Remark~6 in Section~6 for a somewhat 
more detailed statement) on the space of
stationary probability measures
on $(A^\infty\times\hat{A}^\infty)$
equipped with the
topology of weak convergence.
Assume, moreover,
that this LDP holds with a convex,
good rate-function $I(\cdot)$.
[See \cite{dawson-gartner:87}\cite[Sec.~5.3,~5.4]{deuschel-stroock:book}%
\cite[Sec.~6.5.3]{dembo-zeitouni:book}\cite{bryc-dembo:96} 
for a general discussion as well as specific examples 
of processes for which the above conditions
hold. Apart from the $\iid$ case, these
examples also include all ergodic finite-state
Markov chains, among many others.]

It is easy to check that, when 
$\rho$ is bounded and continuous on 
$A \times \hat{A}$, then with
$\BBP$-probability one, conditional on 
$x_1^\infty$, the random variables
$\{\rho_n(x_1^n,Y_1^n)\}$ 
satisfy the LDP upper bound with respect 
to the deterministic, convex rate-function 
$J(D)=\inf  I(\nu)$, where the infimum
is over all stationary probability measures
$\nu$ on $A^\infty \times \hat{A}^\infty$ such that
the $A^\infty$-marginal of $\nu$ is $\BBP$ and
$E_\nu[\rho(X_1,Y_1)] = D$.
Indeed, Comets \cite{comets:89} provides
such an argument when $\Xp$ and $\Yp$ are both $\iid$
Moreover, he shows that in that case 
the corresponding LDP lower bound also holds, 
and hence Theorem 4 applies.
Unfortunately, the conditional
LDP lower bound has to be verified on a 
case-by-case basis.

\medskip

{\em Remark 4:} 
Although quite strong,
the $\psipm$-mixing condition
of Theorem~3, and the $(S)$-mixing
condition of \cite{chi-AP:01},
probably cannot be significantly 
relaxed: For example, in the special case 
when $\Xp$ is a constant process 
taking on just a single value,
if Theorem~3 were to hold (for any bounded 
distortion measure) with a strictly 
monotone rate-function, then necessarily 
the empirical measures of $Y_1^n$ 
would satisfy the LDP in the space 
${\cal P}_a(\Ahat)$ 
(see \cite{bryc-dembo:96} for details).
But
\cite[Example~1]{bryc-dembo:96} illustrates
that this LDP
may fail even when $\Yp$ is a stationary 
ergodic Markov chain with 
discrete alphabet $\Ahat$. In particular,
the example in \cite{bryc-dembo:96}
has an exponential $\phi$-mixing rate. 


\subsection{Generalized AEP for Optimal Lossy Compression}

Here we present a version of the generalized
AEP that is useful in proving direct coding
theorems. Let $\Xp$ be a stationary
ergodic process. For the distortion measure
$\rho$ we adopt two simple regularity conditions.
We assume the existence of a {\em reference 
letter}, i.e., an $\hat{a}\in\Ahat$ such that 
$$E_{P_1}[\rho(X_1,\hat{a})]<\infty.$$
Also, following \cite{kieffer:91}, we 
require that for any distortion level 
$D>0$ there is a scalar quantizer 
for $\Xp$ with finite rate.

\smallskip

{\em Quantization Condition:} 
For each $D>0$
there is a ``quantizer'' $q:A\to B$ for
some countable (finite or infinite)
subset $B\subset\Ahat$, 
such that:

i. $\;\;\rho(x,q(x))\leq D$ for all $x\in A$, and

ii. $\;$
the entropy $H(q(X_1))<\infty$.

\smallskip

\noindent
The following was implicitly proved
in \cite{kieffer:91}; 
see also \cite{konto-zhang:00} 
for details.

\smallskip
 
{\em Theorem~5. Generalized AEP for Optimal Lossy Compression
\cite{kieffer:91}:}
Let $\Xp$ be a stationary ergodic process.
Assume that the distortion measure $\rho$
satisfies the quantization condition,
that a reference letter exists, and 
that for each $n\geq 1$ the infimum of
$$E_{P_n}[-\log Q_n(B(X_1^n,D))]$$
over all probability measures $Q_n$
on $\Ahatn$ is achieved by some 
$\widetilde{Q}_n$.
Then for any $D>0$
\be
-\frac{1}{n}\log \widetilde{Q}_n(B(X_1^n,D)) \to R(D)
        \;\;\;\;\mbox{w.p.1}
\label{eq:chi}
\ee
where $R(D)$ is the rate-distortion 
function of the process $\Xp$.
 
\medskip

{\em Historical Remarks:}
The relevance of the quantities 
$-\log Q_n(B(X_1^n,D))$ to 
information theory was first 
suggested implicitly
by Kieffer \cite{kieffer:91}
and more explicitly 
by {\L}uczak and Szpankowski 
\cite{luczak-szpankowski}.
Since then, many papers have 
appeared proving the generalized AEP 
under different conditions;
we mention here a subset
of those proving some of
the more general results.
The case of finite alphabet
processes was considered by Yang and Kieffer
\cite{yang-kieffer:1}.
The generalized AEP for 
processes with general 
alphabets and $\Yp$ $\iid$ 
was proved by Dembo and
Kontoyiannis
\cite{dembo-kontoyiannis}
and by Yang and Zhang 
\cite{yang-zhang:99}.
Finally, the case when
$\Yp$ is not $\iid$ was
(Theorem~3) treated by
Chi \cite{chi-it:01}\cite{chi-AP:01}.
The observations of Theorem~4 
about the rate-function 
$R(\BBP,\BBQ,D)$ are new.
Theorem~5 essentially
comes from Kieffer's work
\cite{kieffer:91};
see also \cite{konto-zhang:00}.

We should also mention 
that, in 
a somewhat different context,
the intimate relationship
between the AEP and large
deviations is discussed in
some detail by Orey in
\cite{orey:85b}.

\subsection{Densities vs. Balls}
Let us recall the classical generalization
of the AEP, due to Barron \cite{barron:1}
and Orey \cite{orey:85}, to processes with
values in general alphabets. Suppose $\Xp$
as above is a general stationary ergodic process
with marginals $\{P_n\}$ that are
absolutely continuous with respect to
the sequence of measures $\BBM=\{M_n\}$.

\medskip

{\em Theorem~6. AEP for Processes with Densities
\cite{barron:1}\cite{orey:85}:}
Let $\Xp$ be a stationary ergodic process whose
marginals $P_n$ have densities $f_n=dP_n/dM_n$ 
with respect to the $\sigma$-finite measures $M_n$,
$n\geq 1$. Assume that the sequence $\BBM$
of dominating measures is Markov of finite order, 
with a stationary transition measure, and that the 
relative entropies
\ben
H_n\bydef E_{P_n}\left[\log\frac{f_n(X_1^n)}
				  {f_{n-1}(X_1^{n-1})}
		\right],\;\;\;\;n\geq 2,
\een
have $H_n>-\infty$ eventually. Then
\be
-\frac{1}{n}\log \frac{dP_n}{dM_n}(X_1^n)\to -H(\BBP\|\BBM)
        \;\;\;\;\mbox{w.p.1}
\label{eq:BarronAEP}
\ee
where $H(\BBP\|\BBM)$ is the relative entropy
rate defined as $H(\BBP\|\BBM)=\lim_n H_n=\inf_n H_n$.

\medskip

The AEP for processes with densities is
also know to hold when the reference measures
$M_n$ do not form a Markov
sequence, under some additional
mixing conditions (see \cite{orey:85} where
$M_n$ are taken to be non-Markov measures
satisfying an additional mixing condition,
and the more recent extension 
in \cite{chazottesetal:98}
where the $M_n$ are taken to be
discrete Gibbs measures.)
Moreover, Kieffer 
\cite{kieffer:73}\cite{kieffer:73b}
has given counterexamples
illustrating that without some mixing
conditions on $\{M_n\}$ the AEP
(\ref{eq:BarronAEP}) fails to hold.

There is a tempting analogy between
the generalized AEP (\ref{eq:thm4})
and the AEP for processes with 
densities (\ref{eq:BarronAEP}).
The formal similarity between
the two suggests that, if we identify 
the measures $Q_n$ with the reference
measures $M_n$, corresponding results
should hold in the two cases. 
Indeed, this does in general appear
to be the case, as is illustrated 
by the various generalized AEPs
stated above. Moreover, we can
interpret the result of Theorem~5
as the natural analog of the classical 
discrete AEP (\ref{eq:discreteAEP}) to
the case of lossy data compression.
As we argued in the
introduction,
the generalized AEPs of the previous
sections play analogous roles in the proofs
of the corresponding direct coding
theorems.

Taking this analogy further indicates
that there might be a relationship
between these two different 
generalizations.
In particular, when $n$ is large and
the distortion level $D$ is small,
the following heuristic
calculation seems compelling:

\ben
-H(\BBP\|\BBQ)
&\approxa&
-\frac{1}{n}\log \frac{dP_n}{dQ_n}(X_1^n)\\
&\approxb&
-\frac{1}{n}\log \frac{P_n(B(X_1^n,D))}{Q_n(B(X_1^n,D))}\\
&=&
-\frac{1}{n}\log P_n(B(X_1^n,D))
+\frac{1}{n}\log Q_n(B(X_1^n,D))\\
&\approxc& 
R(\BBP,\BBP,D)-R(\BBP,\BBQ,D)\\
&\approxd&-H(\BBP\|\BBQ)
\een
where $(a)$ holds in the limit as $n\to\infty$
by Theorem~6, $(b)$ should hold when $D$ is small
by the assumption that $P_n$ has a density
with respect to $Q_n$, $(c)$ would 
follow in the limit as $n\to\infty$ by 
an application of the generalized AEP, 
and it is natural to conjecture
that $(d)$ holds in the limits
as $D\downarrow 0$ by reading
the above calculation backwards.

In the following two sections we
formalize the above heuristic 
argument in two special cases: 
First when $\Xp$ is a discrete 
process taking values in a finite 
alphabet, and second when $\Xp$ is
a continuous process
taking values in $\RL^d$.

\subsubsection{Discrete Case}
Here we take $\Xp$ to be a stationary ergodic
process taking values in a finite alphabet $A$,
and $\Yp$ to be $\iid$ with first order marginal 
distribution $Q=Q_1$ on the same alphabet $A=\Ahat$.
Similarly we write $P=P_1$ for the first order
marginal of $\Xp$.
In Theorem~7 we justify the above
calculation by showing that the 
limits as $D\downarrow 0$ and as $n\to\infty$
can indeed be taken together in any 
fashion: We show that the double 
limit of the central expression 
\be
r_n(X_1^n,D)
\bydef
\frac{1}{n}\log \frac{P_n(B(X_1^n,D))}{Q_n(B(X_1^n,D))}
\label{eq:ratio}
\ee
is equal to $H(\BBP\|\BBQ)$ with probability 1,
independently of how $n$ grows and 
$D$ decreases to zero. Its proof is 
given in Appendix~A.

\medskip

{\em Theorem~7. Densities vs. Balls in the Discrete Case:}
Let $\Xp$ be a stationary ergodic process
and $\Yp$ be $\iid$, both on the finite 
alphabet $A$. Assume that $\rho(x,y)=0$
if and only if $x=y$, and $Q(x)>0$ for all $x$. 
Then the following
double limit exists:
$$\limnd
\frac{1}{n}\log
        \frac{P_n(B(X_1^n,D))}
             {Q^n(B(X_1^n,D))}
	\;=\; H(\BBP\|\BBQ)
\;\;\;\;\mbox{w.p.1}
$$
In particular, the repeated limit
$\lim_{n}\lim_{D}$
exists with probability one
and is equal to $H(\BBP\|\BBQ)$.

\subsubsection{Continuous Case}

Here we state a weaker version of Theorem~7 in the 
case when
$A=\Ahat=\RL^d$ for some $d\geq 1$, and
when $\Xp$ is an $\RL^d$-valued, 
stationary ergodic process.
Suppose that the marginals $\{P_n\}$ of 
$\Xp$ are absolutely continuous 
with respect to a sequence
of reference measures $\{Q_n\}$. Throughout 
this section we take the $Q_n$
to be product measures, $Q_n=Q^n,$
for some fixed Borel probability
measure $Q$ on $\RL^d$.
A typical example to keep in mind 
is when $Q$ a Gaussian measure on 
$\RL$ and $\Xp$ a real-valued stationary 
ergodic process all of whose marginals 
$P_n$ have continuous densities
with respect to Lebesgue measure.

For simplicity, we take $\rho$ to be 
squared-error distortion,
$\rho(x,y)=(x-y)^2$, although
the proof of Theorem~8, given in 
Appendix~B, may easily be adapted to 
apply for somewhat more general
difference distortion measures.

\medskip

{\em Theorem~8. Densities vs. Balls in the Continuous Case:}
Let $\Xp$ be an $\RL^d$-valued stationary ergodic process,
whose marginals $P_n$ have densities $f_n=dP_n/dQ_n$ with 
respect to a sequence of product measures $Q_n=Q^n$, 
$n\geq 1$, for a given probability measure $Q$ on $\RL^d$.
Let $\rho(x,y)=(x-y)^2$ for any $x,y \in \RL^d$.

(a) The following repeated limit holds:
$$\lim_{n\to\infty}
  \lim_{D\downarrow 0}\;
        \frac{1}{n}\log
        \frac{P_n(B(X_1^n,D))}
             {Q_n(B(X_1^n,D))}
	= H(\BBP\|\BBQ)
	\;\;\;\;\mbox{w.p.1.}
$$

(b) Assume, moreover, that $\Xp$ is $\iid$ 
with marginal distribution $P_1=P$ on $\RL^d$, 
and that the following conditions are satisfied:
Both $E_{P\times Q}[\rho(X,Y)]$ and
$E_{P\times P}[\rho(X,Y)]$ are finite 
and nonzero; the expectation
$$E_P[-\log Q(B(X,D))] 
\;\;\;\;\mbox{is finite for all}\;D>0;$$
and a $\delta>0$ exists for which
\be
E_P\left[\sup_{0<D<\delta} \left|
	\log \frac{P(B(X,D))}{Q(B(X,D))} \right|\right]<\infty.
\label{eq:integrability}
\ee
Then, the reverse repeated limit also holds:
$$\lim_{D\downarrow 0}
  \lim_{n\to\infty}\;
        \frac{1}{n}\log
        \frac{P_n(B(X_1^n,D))}
             {Q_n(B(X_1^n,D))}
	= H(\BBP\|\BBQ)
	\;\;\;\;\mbox{w.p.1.}
$$

\medskip

It is easy to check that all conditions of the
theorem hold when $Q$ is a Gaussian measure on $\RL$
and $P$ has finite variance and a probability density 
function $g$ (with respect to Lebesgue measure)
such that $E_P(\sup_{|y-X|<\delta} |\log g(y)|)<\infty$ 
for some $\delta>0$.  For example, this is the case 
when both $P$ and $Q$ are Gaussian distributions on $\RL$.

As will be seen from the proof of the theorem,
although we are primarily interested in the 
case when the relative entropy rate $H(\BBP\|\BBQ)$
is finite, the result remains true when 
$H(\BBP\|\BBQ)=\infty$, and in that case
assumption (\ref{eq:integrability}) can be relaxed to
$$E_P\left[\sup_{0<D<\delta} \log \frac{Q(B(X,D))}{P(B(X,D))} 
\right]<\infty.$$


Finally we note that, in the context of ergodic
theory, Feldman \cite{feldman:80} developed
a different verison of the generalized AEP,
and also discussed the relationship between 
the two types of asymptotics (as $n\to\infty$, 
and as $D\downarrow 0$).

\section{Applications of the Generalized AEP}
As outlined in the
introduction, the generalized AEP can
be applied to a number of problems in data compression
and pattern matching. Following along the lines of the
corresponding applications in the lossless case, below
we present applications of the results of the previous
section to: 1.~Shannon's random coding schemes;
2.~mismatched codebooks in lossy data compression;
3.~waiting times between stationary processes
(corresponding to idealized Lempel-Ziv coding);
4.~practical lossy Lempel-Ziv coding for memoryless
sources; and 5.~weighted codebooks in 
rate-distortion theory.
 
\subsection{Shannon's Random Codes}
Shannon's well-known construction of optimal
codes for lossy data compression is based on
the idea of generating a random codebook. We
review here a slightly modified version of
his construction \cite{shannon:59}
and describe how the performance of the
resulting random code can be analyzed 
using the generalized AEP.

Given a sequence of probability distributions
$Q_n$ on $\Ahatn$, $n\geq 1$, we generate a 
{\em random codebook according to the measures $Q_n$}
as an infinite sequence of $\iid$ random vectors
$$Y_1^n(i),\;\;\;\;i\geq 1$$
with each $Y_1^n(i)$ having distribution
$Q_n$ on $\Ahatn$. Suppose that, for a fixed $n$,
this codebook is available to both the encoder and decoder. 
Given a source string $X_1^n$ to
be described with distortion $D$ or less,
the encoder looks for a $D$-close match of
$X_1^n$ into the codebook $\{Y_1^n(i)\;;\;i\geq 1\}$.
Let $i_n$ be the position of the first such match
\ben
i_n\bydef \inf \{i\geq 1\;:\;\rho_n(X_1^n,Y_1^n(i))\leq D\}
\een
with the convention that the infimum of
the empty set equals $+\infty$. If a
match is found, then the encoder describes
to the decoder the position $i_n$ using
Elias' code for the integers 
\cite{elias}. This takes no more than
\be
\log_2 i_n + 2\log_2\log_2 i_n + \mbox{Const.}
\;\;\;\;\mbox{bits}.
\label{eq:elias}
\ee
If no match is found 
(something that asymptotically will 
{\em not} happen, with probability one),
then the encoder describes $X_1^n$ with
distortion $D$ or less using some other
default scheme. 

Let $\ell_n(X_1^n)$
denote the overall description 
length of the algorithm just 
described. In view of (\ref{eq:elias}),
in order to understand its 
compression performance,
that is, to understand the 
asymptotic behavior of 
$\ell_n(X_1^n)$, it suffices 
to understand the behavior of the quantity
$$\log_2 i_n,\;\;\;\;\mbox{for large $n$.}$$
Suppose that the probability
$Q_n(B(X_1^n,D))$ of finding a $D$-close
match for $X_1^n$ in the codebook is nonzero.
Then, conditional on the source string $X_1^n$, 
the distribution of $i_n$ is geometric with
parameter $Q_n(B(X_1^n,D))$. From this 
observation is easy to deduce that 
the behavior of $i_n$ is closely
related to the behavior of the quantity
$1/Q_n(B(X_1^n,D))$. The next theorem is
an easy consequence of this fact so it is
stated here without proof; see the 
corresponding arguments in 
\cite{kontoyiannis-red:00}\cite{konto-zhang:00}.

\medskip

{\em Theorem~9. Strong Approximation:}
Let $\Xp$ be an arbitrary process and
let $\{Q_n\}$ be a given sequence of 
codebook distributions. 
If $Q_n(B(X_1^n,D))>0$ eventually with 
probability one,
then for any $\epsilon>0$:
\ben
\log_2 i_n 
&\leq& -\log_2 Q_n(B(X_1^n,D)) + \log_2\log_2 n + 3
	\;\;\;\;\mbox{eventually, w.p.1}\\
\mbox{and}\;\;
	\log_2 i_n 
&\geq&
-\log_2 Q_n(B(X_1^n,D)) -
\log_2 n - (1+\epsilon)\log_2\log_2 n
\;\;\;\;\mbox{eventually, w.p.1.}
\een

\medskip

The above estimates can now be combined 
with the results of the generalized AEP 
in the previous section to determine the
performance of codes based on random 
codebooks with respect to the ``optimal''
measures $Q_n$. To illustrate this 
approach we consider the special case
of memoryless sources and finite
reproduction alphabets, and show that
the random code with respect to 
(almost) any random codebook realization 
is asymptotically optimal, with 
probability one. Note that corresponding
results can be proved, in exactly the
same way, under much more general 
assumptions. For example, utilizing
Theorem~5 instead of Theorem~1 we
can prove the analog of Theorem~10
below for arbitrary stationary 
ergodic sources.

Let $\Xp$ be an $\iid$ source with 
marginal distribution $P_1=P$ 
on $A$, and take the reproduction
alphabet $\Ahat$ to be finite.
For simplicity we will
assume that the distortion measure 
$\rho$ is bounded, i.e., 
$\sup_{x,y}\rho(x,y)<\infty,$
and we also make the customary 
assumption that 
\be
\sup_{x\in A}\min_{y\in\hat{A}}\rho(x,y) = 0.
\label{eq:maximin}
\ee
[See the remark at the end of Section~5.1.1 for
a discussion of this condition and when it can 
be relaxed.]
As usual, we define the rate-distortion
function of the memoryless source
$\Xp$ by
$$R(D)=\inf_{(X,Y)}\,I(X;Y)$$
where the infimum is over all jointly 
distributed random variables $(X,Y)$ 
with values in $A\times\Ahat$, such 
that $X$ has distribution $P$
and $E[\rho(X,Y)]\leq D$.
Let 
\be
\Dbar\bydef \min_{y\in\hat{A}}E_P[\rho(X,y)]
\label{eq:Dbar}
\ee
and note that $R(D)=0$ for $D\geq\Dbar$.
To avoid the trivial case when
$R(D)=0$ for all $D,$ we assume
that $\Dbar>0$ and we restrict 
our attention to the interesting
range of values $D\in(0,\Dbar)$.
Recall \cite{yang-zhang:99}\cite{kontoyiannis-red:00}
that for any such $D$,
$R(D)$ can alternatively be 
written as
$$R(D)=\inf_Q R_1(P,Q,D)$$
where the infimum is over all 
probability distributions $Q$ on $\Ahat$.
Since we take $\Ahat$ to be finite,
this infimum is always achieved
(see \cite{kontoyiannis-red:00})
by a probability distribution
$Q=Q^*$. 
To avoid cumbersome notation in the statements
of the coding theorems given next and also in 
later parts of the paper, we also write
$\calR(D)$ for the rate-distortion 
function of the source $\Xp$ expressed
in {\em bits} rather than in nats:
$$\calR(D)\bydef (\log_2 e)R(D).$$
Finally, we write $Q_n^*$ for the product
measures $(Q^*)^n$ and call
$\{Q_n^*\}$ the {\em optimal reproduction 
distributions at distortion level $D$.}

Combining Theorem~9 with the 
generalized AEP of Theorem~1 
implies the following 
strengthened direct
coding theorem.

\medskip
 
{\em Theorem~10. Pointwise Coding Theorem 
for I.I.D. Sources \cite{kontoyiannis-red:00}:}
Let $\Xp$ be an $\iid$ source with distribution
$P$ on $A$, and let $Q_n^*$ denote the optimal
reproduction distributions at distortion level
$D\in(0,\Dbar)$.
Then the codes based on almost any realization 
of the Shannon random codebooks according
to the measures $\{Q_n^*\}$ have codelengths 
$\ell_n(X_1^n)$
satisfying:
$$\lim_{n\to\infty}\frac{1}{n}\ell_n(X_1^n)
= \calR(D)\;\;\;\;\mbox{bits per symbol, w.p.1.}$$

\medskip

A simple modification of the above scheme can
be used to obtain {\em universal} codebooks
that achieve optimal compression for any 
memoryless source:
Given a fixed block-length $n$, we consider
the collection of all $n$-types on $\Ahat$,
namely, all distributions $Q$ of the form
$Q(\hat{a})=j/n$, $0\leq j\leq n$, for 
$\hat{a}\in\Ahat$. Instead of generating 
a single random codebook according to the 
optimal distribution $Q_n^*$, we generate 
{\em multiple codebooks}, one for each
product measure $Q^n$ corresponding to an
$n$-type $Q$ on $\Ahat$. Then we (as the 
encoder) adopt a greedy coding strategy. We find 
the first $D$-close match for $X_1^n$ in 
each of the codebooks, and pick the one 
in which the match appears the earliest.
To describe $X_1^n$ to the decoder with
distortion $D$ or less we then describe
two things: (a)~the index of the codebook 
in which the earliest match was found, 
and (b)~the position $i_n$ of this 
earliest match. Since there are at
most polynomially many $n$-types
(cf. \cite{csiszar:book}\cite{cover:book}),
the rate of the description of (a) is 
asymptotically negligible. Moreover,
since the set of $n$-types is 
asymptotically dense among probability
measures on $\Ahat$, we eventually
do as well as if we were using the 
optimum codebook distribution $Q_n^*$.

\medskip
 
{\em Theorem~11. Pointwise Universal Coding Theorem
\cite{kontoyiannis-red:00}:}
Let $\Xp$ be an arbitrary $\iid$ source with 
distribution $P$ on $A$, let $R(D)$
be the rate-distortion function of this source
at distortion level $D\in(0,\Dbar)$,
and let $\calR(D)$ denote its 
rate-distortion function in bits.
The codes 
based on almost any realization of the 
universal Shannon random codebooks have 
codelengths $\ell_n(X_1^n)$ satisfying:
$$\lim_{n\to\infty}\frac{1}{n}\ell_n(X_1^n)
= \calR(D)\;\;\;\;\mbox{bits per symbol, w.p.1.}$$

\subsection{Mismatched Codebooks}
In the last section we described how,
for memoryless sources, the Shannon 
random codebooks with respect to the
optimal reproduction distributions can
be used to achieve asymptotically 
optimal compression performance. In this 
section we briefly consider the question
of determining the rate achieved
when an arbitrary (stationary ergodic)
source $\Xp$ is encoded using a 
random codebook according to the 
$\iid$ distributions $Q^n$, 
for an arbitrary distribution $Q$ 
on $\Ahat$. For further discussion 
of the problem of mismatched
codebooks see 
\cite{sakrison:69}\cite{sakrison:70}\cite{lapidoth:97}\cite{kanlis:phd}
and the references therein.

The following theorem is an immediate
consequence of combining Theorem~1  
with Theorem~9 and the discussion in 
Section~3.1 (see also Example~1 in 
Section~2.2).

\medskip

{\em Theorem~12. Mismatched Coding Rate:}
Let $\Xp$ be a stationary ergodic process
with marginal distribution $P_1=P$ on $A$,
let $Q$ be an arbitrary distribution
on $\Ahat$, and define $\Dmin$ and
$\Dav$ as in Section~2.2.
\begin{itemize}
\item[(a)]{\em Arbitrary I.I.D. Codebooks:}
  For any distortion level $D\in(\Dmin,\Dav)$,
  the codes based on almost any realization
  of the Shannon random codebooks according
  to the measures $\{Q^n\}$ have codelengths 
	$\ell_n(X_1^n)$ satisfying:
  $$\lim_{n\to\infty}\frac{1}{n}\ell_n(X_1^n)
  = (\log_2 e)R_1(P,Q,D)\;\;\;\;\mbox{bits per symbol, w.p.1.}$$
\item[(b)]{\em I.I.D. Gaussian Codebooks:}
  Suppose
$\rho(x,y)=(x-y)^2$ and
$\Xp$ is a real-valued process with
  finite variance $\sigma^2=\VAR(X_1)$.
  Let $Q$ be the $N(0,\tau^2)$ distribution 
  on $\RL$. Then for any distortion level 
  $D\in(0,\sigma^2+\tau^2)$, the codes based on 
  almost any realization of the Gaussian
  codebooks according to the measures $\{Q^n\}$ 
  have codelengths 
  $\ell_n(X_1^n)$ satisfying:
  $$\lim_{n\to\infty}\frac{1}{n}\ell_n(X_1^n)
  = \frac{1}{2}\log_2\left(\frac{v}{D}\right)
	-(\log_2 e)\frac{(v-D)(v-\sigma^2)}
                        {2v\tau^2}
	\;\;\;\;\mbox{bits per symbol, w.p.1,}$$
  where
  $$v\bydef\frac{1}{2}\left[\tau^2+\sqrt{\tau^4+4D\sigma^2}\right].$$
\end{itemize}

\medskip

{\em Lossless vs. Lossy Mismatch:} Recall
that, in the case of lossless data compression,
if instead of the true source distribution
$P$ a different coding distribution $Q$ is used,
then the code-rate achieved is 
\be
H(P)+H(P\|Q).
\label{eq:penalty1}
\ee
Similarly in the current setting of lossy 
data compression, if instead of the optimal
reproduction distribution $Q^*$ we use a
different codebook distribution $Q$, the
rate we achieve is $R_1(P,Q,D)$. 
An upper bound for $R_1(P,Q,D)$ is
obtained by taking $(X,Y)$
in the expression of Remark~1 
to be the jointly distributed random 
variables that achieve the infimum 
in the definition of the rate-distortion
function of $P$. Then the (mismatched) rate
of the random code based on $Q$ instead 
of $Q^*$ is:
\be
R_1(P,Q,D)\leq R(D) + H(Q^*\|Q).
\label{eq:penalty2}
\ee
Equations (\ref{eq:penalty1})
and (\ref{eq:penalty2}) illustrate
the analogy between the penalty terms 
in the lossless and lossy case
due to mismatch.

\medskip

Next we discuss two special cases 
of part~(b) of the theorem
that are of particular interest.

\medskip

{\em Example~2: Gaussian codebook with 
mismatched distribution:}
Consider the following coding scenario:
We want to encode data generated by 
an $\iid$ Gaussian process 
with $N(0,\sigma^2)$ distribution,
with squared-error distortion
$D$ or less.
In this case, it is well-known 
\cite{berger:book}\cite{cover:book} that
for any $D\in(0,\sigma^2)$ 
the optimal reproduction distribution $Q^*$
is the $N(0,\sigma^2-D)$ distribution,
so we construct random codebooks 
according to the $\iid$ distributions
$Q_n^*=(Q^*)^n$.

But suppose that, instead of an
$\iid$ Gaussian, the source turns out
to be some arbitrary stationary ergodic 
$\Xp$ with zero mean and variance $\sigma^2$.
Theorem~12~(b) implies that the asymptotic
rate achieved by our $\iid$ Gaussian
codebook is equal to
$$\frac{1}{2}\log_2\left(\frac{\sigma^2}{D}\right)
\;\;\;\;\mbox{bits per symbol.}$$
Since this is exactly the
rate-distortion function of the 
$\iid$ $N(0,\sigma^2)$ source, we
conclude that the rate achieved is
the same as what we would have 
obtained on the Gaussian source we
originally expected. This offers
yet another justification of the
folk theorem that the Gaussian 
source is the hardest one to compress,
among sources with a fixed variance. 
In fact, the above result is 
a natural fixed-distortion
analog of \cite[Theorem~3]{lapidoth:97}.

\medskip

{\em Example~3: Gaussian codebook with mismatched variance:}
Here we consider a different type of mismatch.
As before, we are prepared to encode an 
$\iid$ Gaussian source, but we have an 
incorrect estimate of its variance,
say $\hat{\sigma}^2$ instead of the true
variance $\sigma^2$. So we are using
a random codebook with respect to the
optimal reproduction distribution
$Q_n^*=(Q^*)^n$, where $Q^*$ is the
$N(0,\hat{\sigma}^2-D)$ distribution,
but the actual source is $\iid$ 
$N(0,\sigma^2)$. In this case, 
the rate achieved by 
the random codebooks according to 
the distributions $Q_n^*$ is given
by the expression in Theorem~12~(b),
with $\tau^2$ replaced by $\hat{\sigma}^2-D$.
Although the resulting expression 
is somewhat long and not easy to 
manipulate analytically, it is 
straightforward to evaluate 
numerically. For example, 
Figure~1 shows the asymptotic 
rate achieved, as a function of the 
error $e=\sigma^2-\hat{\sigma}^2$
in the estimate of the true variance.
As expected, the best rate is 
achieved when the codebook distribution 
is matched the source (corresponding to $e=0$),
and it is equal to the rate-distortion function
of the source. Moreover, as one might 
expect, it is more harmful to 
underestimate the variance 
than to overestimate it.

\begin{figure}[ht]
\centerline{\epsfxsize 3.2in \epsfbox{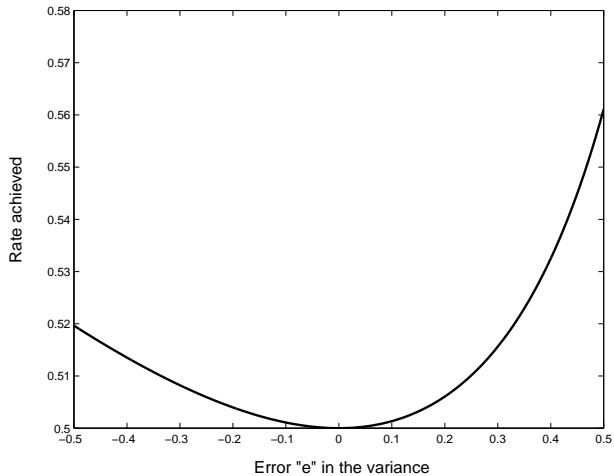}}
\caption{This graph shows the rate achieved by an
$\iid$ Gaussian codebook of variance $\hat{\sigma}^2-D$
when applied to $\iid$ $N(0,\sigma^2)$ data.
The rate is shown as a function of the error
$e=\sigma^2-\hat{\sigma}^2$ in the variance estimate.
In this particular example: $\sigma^2=2$, $D=1$, 
the error $e$ ranges from $-1/2$ to $1/2$,
and the rate-distortion function of the source
equals 0.5 bits/symbol.}
\end{figure}

\subsection{Waiting Times and Idealized Lempel-Ziv Coding}
Given $D\geq 0$ and two independent realizations
from the stationary ergodic processes $\Xp$ and $\Yp$, 
our main quantity of interest here is the 
{\em waiting time} $W_n=W_n(D)$ until a $D$-close 
version of the initial string $X_1^n$ first appears
in $Y_1^\infty$. Formally
\be
W_n\;=\;\inf\{i\geq 1\; :\; 
	\rho_n(X_1^n,Y_i^{i+n-1})\leq D\}
\label{eq:Wn-def}
\ee
with the convention, as before, that the infimum 
of the empty set equals $+\infty$.

The motivation for studying the asymptotic behavior
of $W_n$ for large $n$ is twofold.

\medskip

{\em Idealized Lempel-Ziv coding.}
The natural extension of the idealized
scenario described in the
introduction
is to consider a message $X_1^n$ that
is to be encoded with the help of a
database $Y_1^\infty$.
The
source and
the database are assumed to be 
independent, and the database
distribution may or may not be
the same as that of the source.
In order to communicate $X_1^n$
to the decoder with distortion
$D$ or less, the encoder simply
describes $W_n$, using no more
than
$$\log_2 W_n + O(\log_2\log_2 W_n)\;\;\;\;\mbox{bits.}$$
Therefore, the asymptotic performance 
of this idealized scheme can be 
completely understood in terms
of the asymptotics of $\log W_n$,
for large $n$.

\medskip

{\em DNA pattern matching.}
Here we imagine that $X_1^n$ represents 
a DNA or protein ``template,'' and we 
want to see whether it appears, either
exactly or approximately, as a contiguous 
substring of a database DNA sequence 
$Y_1^\infty$. We are interested in
quantifying the ``degree of surprise''
in the fact that a $D$-close match was
found at position $W_n$. Specifically,
was the match found ``atypically''
early, or is the value of $W_n$ 
consistent with the hypothesis
that the template and the database
are independent? For a detailed
discussion, see, e.g.,
\cite[Section~3.2]{dembo-zeitouni:book}\cite{karlin-ost:88}%
\cite{agw:90}\cite{arratia-waterman}
and the references therein.


\newpage

If for a moment we consider
the case when both $\Xp$ and
$\Yp$ are $\iid$, we see that 
the waiting time $W_n$ is,
at least intuitively, closely
related to the index $i_n$ of 
Section~3.1.
As the following result shows,
although the distribution of
$W_n$ is not exactly geometric,
$W_n$ behaves very much 
like $i_n$, at least in the 
exponent. That is, the
difference
$$\log W_n -[-\log Q_n(B(X_1^n,D))]$$
is ``small,'' eventually
with probability one.

Recall the definition of
$\psi$-mixing from Section~2.3, and
also the definition of the
$\phi$-mixing coefficients of $\Yp$
$$\phi(k)\;=\;\sup\{|\BBQ(B|A)-\BBQ(B)|\;:\;\;
B\in\sigma(Y_{k}^{\infty}),\;
A\in\sigma(Y_{-\infty}^0),\; \BBQ(A)>0\}$$
where, as before, $\sigma(Y_i^j)$ denotes
the $\sigma$-field generated by $Y_i^j$. 
The process $\Yp$ is called
{\em $\phi$-mixing}
if $\phi(k)\to 0$ as $k\to\infty$;
see \cite{bradley} for an extensive 
discussion of $\phi$-mixing and related
mixing conditions.
 
\medskip

{\em Theorem~13. Strong Approximation
\cite{kontoyiannis-jtp}\cite{dembo-kontoyiannis}:}
Let $\Xp$ and $\Yp$ be stationary ergodic processes,
and assume that $\Yp$ is either $\psi$-mixing
or $\phi$-mixing with summable $\phi$-mixing
coefficients, $\sum_{k\geq 1} \phi(k)<\infty$.
If $Q_n(B(X_1^n,D))>0$ eventually with 
probability one,
then for any $\epsilon>0$:
\ben
-(1+\epsilon)\log n
\;\leq\;
\log [W_n Q_n(B(X_1^n,D))]
\;\leq\;
(2+\epsilon)\log n
\;\;\;\;\mbox{eventually, w.p.1.}
\een

\medskip

Theorem~13 of course implies that
\be
\log W_n = -\log Q_n(B(X_1^n,D)) + O(\log n)
\;\;\;\;\mbox{w.p.1}
\label{eq:strong}
\ee
and combining this with the generalized
AEP statements of Theorems~1 and~4 we
immediately obtain the first order 
(or strong-law-of-large-numbers, SLLN)
asymptotic behavior of the waiting
times $W_n$:

\medskip

{\em Theorem~14. SLLN for Waiting Times:}
Let $\Xp$ and $\Yp$ be stationary ergodic processes.

(a)~If $\Yp$ is $\iid$ and the 
average distortion $\Dav$ is finite,
then for any $D\in(\Dmin,\Dav)$
\be
\frac{1}{n}\log W_n \to R_1(P_1,Q_1,D)
\;\;\;\;\mbox{w.p.1.}
\label{eq:w-slln}
\ee

(b)~If $\Yp$ is $\psi$-mixing and the distortion
	measure $\rho$ is bounded, then for any
	$D\in(\Dmin,\Dav)$
\be
\frac{1}{n}\log W_n \to R(\BBP,\BBQ,D)
\;\;\;\;\mbox{w.p.1.}
\label{eq:w-slln2}
\ee

\medskip

Note that similar results can be
obtained under different assumptions
on the process $\Yp$, using Theorems~3 
and~5 in place of Theorems~1 and~4 as
done above.
When $\Xp$ is taken to be an
arbitrary stationary ergodic process,
it is natural to expect that the
mixing conditions for $\Yp$ in
Theorem~14~(b) cannot be 
substantially relaxed.
In fact, even in the case of exact matching 
between finite-alphabet processes, Shields 
\cite{shields:3}
has produced a counterexample demonstrating
that the analog of Theorem~13 does not hold 
for arbitrary stationary ergodic $\Yp$. 

\medskip

{\em Historical Remarks:}
Waiting times in the context
of lossy data compression were
studied by Steinberg and Gutman 
\cite{steinberg-gutman} and {\L}uczak 
and Szpankowski \cite{luczak-szpankowski}.
Yang and Kieffer \cite{yang-kieffer:1}
identified the limiting rate-function
for a wide range of finite alphabet
sources, and Dembo and Kontoyiannis 
\cite{dembo-kontoyiannis} and Chi 
\cite{chi-it:01}
generalized these results to processes with
general alphabets.

The strong approximation idea was 
introduced
in \cite{kontoyiannis-jtp}
in the case of exact matching.  For
processes $\Yp$ with summable
$\phi$-mixing coefficients, Theorem~13 was 
proved in \cite{dembo-kontoyiannis}, and when 
$\Yp$ is $\psi$-mixing it was proved, for the 
case of no distortion, in \cite{kontoyiannis-jtp}.
Examining the latter proof, \cite{chi-it:01}
observed that it immediately generalizes to
the statement of Theorem~13.

Related results were obtained by
Kanaya and Muramatsu \cite{kanaya-muramatsu:97},
who extended some of the results of
\cite{steinberg-gutman}
to processes with general alphabets,
and by Koga and Arimoto \cite{koga-arimoto:98}
who considered {\em non-overlapping} waiting 
times between finite-alphabet processes
and Gaussian processes.
Finally, Shields \cite{shields:3}
and Marton and Shields 
\cite{marton-shields:1}
considered waiting times with
respect to Hamming distortion
and for $\Xp$ and $\Yp$ having
the same distribution over a
finite alphabet. For
the case of small
distortion they showed,
under some conditions, 
that approximate matching 
results like (\ref{eq:w-slln}) 
and (\ref{eq:w-slln2}) 
reduce to their natural 
exact matching analogs as
$D\to 0$.

\subsection{Match-Lengths and Practical Lempel-Ziv Coding}
In the idealized coding scenario of
the previous section we considered the
case where a fixed-length message $X_1^n$
is to be compressed using an infinitely 
long database $Y_1^\infty$. But, in practice, 
the reverse situation is much more common:
We typically have a ``long'' message 
$(X_1,X_2,\dots)$ to be compressed, and 
only a finite-length database $Y_1^m$
is available to the encoder and decoder.
It is therefore natural (following
the corresponding development in the
case of lossless compression)
to try and match ``as much as possible'' 
from the message $(X_1,X_2,\dots)$ into the 
database $Y_1^m$. 
With this in mind we 
define the {\em match-length}
$L_m$ as the length $\ell$ of
the longest prefix $X_1^\ell$ that
matches somewhere in the database
with distortion $D$ or less:
\be
L_m=\sup\{\ell \geq 1\;:\;
\rho_\ell(X_1^\ell,Y_{j}^{j+\ell-1})\leq D,
\;\;\mbox{for some}\;\;j=1,2,\ldots,m\}.
\label{eq:Lm-def}
\ee

Intuitively, there is a connection 
between match-lengths and waiting times. 
Long matches should mean short waiting times, 
and vice versa. In the case of exact matching
this connection was precisely formalized by
Wyner and Ziv \cite{wyner-ziv:1}, who observed 
that the following ``duality'' relationship 
always holds:
\be
W_n\leq m
\;\;\;\;
\Leftrightarrow
\;\;\;\;
L_m\geq n.
\label{eq:easy-dual}
\ee
This is almost identical to the
standard relationship in renewal
theory between the number of
events by a certain time and
the time of the $n$th event
(see, e.g., \cite{fellerII:book}).
Wyner and Ziv \cite{wyner-ziv:1}
utilized (\ref{eq:easy-dual})
to translate their first order
asymptotic results about $W_n$ 
to corresponding results about
$L_m$. 

Unfortunately this simple relationship 
no longer holds in the case of
{\em approximate} matching,
when a distortion measure
is introduced. Instead, the following 
modified duality was employed
in \cite{dembo-kontoyiannis} 
to obtain corresponding 
results in approximate matching 
and lossy data compression:
\be
W_n\leq m\;\;\Rightarrow\;\;L_m\geq n
\;\;\;\;
\mbox{and}
\;\;\;\;
L_m\geq n\;\;\Rightarrow\;\;\inf_{k\geq n} W_k\leq m.
\label{eq:duality}
\ee
In \cite{dembo-kontoyiannis} it is shown
that (\ref{eq:duality}) can be used to
deduce the asymptotic behavior of $L_m$
from that of $W_n$,
but this translation
is not straightforward anymore.
In fact, as we discuss in Section~5.2, 
a somewhat more
delicate analysis is needed
in this case.
Nevertheless,
once the
behavior of the waiting
times is understood, 
the first implication in 
(\ref{eq:duality}) immediately 
yields asymptotic {\em lower bounds} 
on the behavior of the match-lengths.
This is significant for data compression
since long match-lengths usually mean
good compression performance.
Indeed, this observation allowed
\cite{kontoyiannis-lossy1-1} to introduce
a new lossy version of the Lempel-Ziv algorithm
that achieves asymptotically optimal 
compression performance for
memoryless sources.
The key characteristics of the 
algorithm are that it has
polynomial implementation 
complexity, and that it
achieves redundancy comparable
to that of its lossless counterpart, 
the FDLZ \cite{wyner-ziv:3}.

We also
mention that, before
\cite{kontoyiannis-lossy1-1},
several practical (yet suboptimal)
lossy versions of the Lempel-Ziv
algorithm were introduced,
perhaps most notably 
by Steinberg and Gutman
\cite{steinberg-gutman} and {\L}uczak 
and Szpankowski \cite{luczak-szpankowski}.
Roughly speaking, the reason for
their suboptimal compression performance 
was that the coding was done with respect
to a database that had the same 
distribution as the source. In view
of the discussion in the previous
section, it is clear that the asymptotic
code-rate of these algorithms is
$R_1(P,P,D)$, which is typically
significantly larger than
the optimal $R(D)=\inf_Q R_1(P,Q,D)$;
see
\cite{yang-kieffer:1} or 
\cite{kontoyiannis-lossy1-1}
for
more detailed discussions.

\subsection{Weighted Codebooks and Sphere-Covering}
Here we describe a related question that was recently 
considered in \cite{covering-TR:99}. In the classical 
rate-distortion problem, one is interested in finding 
``efficient'' codebooks for describing the output of some 
random source to within some tolerable distortion 
level. In terms of data compression, a codebook is
``efficient'' when it contains relatively few codewords,
so that it yields a code with a low rate. Here we are 
interested in the more general problem of finding 
codebooks with small ``mass.'' 

Let $\Xp$ be an
$\iid$ process with marginal distribution
$P$ on a finite alphabet $A$, 
and take $\Ahat = A$ and $\rho$
a distortion measure with the property
that $\rho(x,y)=0$ if and only if $x=y$.
Let $M:A\to(0,\infty)$ be
an arbitrary nonnegative function
assigning mass $M^n(C_n)$ to
subsets $C_n$ of $A^n$:
$$M^n(C_n)\bydef\sum_{y_1^n\in C_n} M^n(y_1^n)
	\bydef\sum_{y_1^n\in C_n}\prod_{i=1}^n M(y_i).$$

The question of interest here can be
stated as follows. Let $C_n$ be 
a subset $A^n$ (we think of $C_n$ 
as the codebook) that
nearly $D$-covers all of $A^n$, 
i.e., with high probability, 
every string $X_1^n$ generated
by the source will match at 
least one element of $C_n$ 
with distortion $D$ or less:
\be
P^n\{\mbox{there is an $y_1^n\in C_n$ such that}\;
	\rho_n(X_1^n,y_1^n)\leq D\}\approx 1.
\label{eq:cover}
\ee
If (\ref{eq:cover}) holds,
how small can the mass of $C_n$ be?
	
For example,
taking $M$ identically equal to one, 
this problem reduces to the rate-distortion
question. Taking $M$ to be a different
probability measure $Q$, it reduces to
the classical hypothesis testing question,
whereas $M=P$ (the source distribution) 
yields ``converses''
to some measure-concentration inequalities;
see \cite{covering-TR:99}
for a detailed treatment of 
these and more general cases.

The next result characterizes the best growth
exponent for the mass of an
arbitrary codebook $C_n$.

\medskip

{\em Theorem~15: Weighted Codebooks \cite{covering-TR:99}:}
Let $\Xp$ be an $\iid$ source on the finite
alphabet $A=\Ahat$, and suppose that 
$\rho(x,y)=0$ if and only if $x=y.$
\begin{itemize}
\item[$(\Leftarrow)$] Let $C_n$ be an arbitrary subset of $A^n$,
and write $D$ for  the expected distance of a source string 
$X_1^n$ from $C_n$:
$$D=E_{P^n}[\min_{y_1^n\in C_n}
	\rho_n(X_1^n,y_1^n)].$$
Then
$$M^n(C_n)\geq e^{nr(D)}$$
where the rate-function $r(D)=r(D;P,M)$ is defined by
$$r(D)=r(D;P,M)=\inf_{(X,Y)}\{I(X;Y)+ E[\log M(Y)]\}$$
and the infimum is taken over all jointly distributed
random variables $(X,Y)$ with values in $A$, such that
$X\sim P$ and $E[\rho(X,Y)]\leq D.$
\item[$(\Rightarrow)$] 
For every $D\geq 0$ 
there is a sequence
of codebooks $\{C^*_n\}$ such that 
\ben
&&
        \limsup_{n\to\infty}\;
        \frac{1}{n}\log M^n(C^*_n)\leq r(D)\\
\mbox{and}&&
        \limsup_{n\to\infty}\;
        E_{P^n}[\min_{y_1^n\in C^*_n}
	\rho_n(X_1^n,y_1^n)]\leq D.
\een
\end{itemize} 

\medskip

The main ingredient in the proof of the direct 
coding theorem in part~$(\Rightarrow)$ above 
is provided by yet another version
of the generalized AEP. Let $(X^*,Y^*)$ be a pair
of random variables achieving the infimum in the
definition of $r(D)$, and let $Q^*$ be the 
distribution of $Y^*$. Now for $\delta>0$ 
and $n\geq 1$ define the sets
$${\cal G}_n=\{y_1^n\in A^n\;:\;
        \hat{P}_{y_1^n}(b)\leq Q^*(b)+\delta,
        \;\;\forall\, b\in A\}$$
where $\hat{P}_{y_1^n}$ denotes the empirical
distribution induced by $y_1^n$ on $A$.
For each $n\geq 1$ define the ``conditioned''
measure $Q^{(c)}_n$ on $A^n$ by conditioning
the product measure $(Q^*)^n$ to the set
${\cal G}_n$. The next theorem provides
the necessary version of the generalized 
AEP in this case.

\medskip

{\em Theorem~16: Generalized AEP 
for Conditioned Measures \cite{covering-TR:99}:}
With the conditioned measures $Q^{(c)}_n$ defined
as above, we have:
$$\limsup_{n\to\infty} -\frac{1}{n}\log Q^{(c)}_n(B(X_1^n,D))
	\leq I(X^*;Y^*) \;\;\;\;\mbox{w.p.1.}$$

\section{Refinements of the Generalized AEP}
As we saw in Section~3, the generalized AEP can
be used to determine the first order asymptotic
behavior of a number of interesting objects 
arising in applications. For example, the 
generalized AEP of Theorem~1
$$-\frac{1}{n}\log Q^n(B(X_1^n,D))\to R_1(P,Q,D)
	 \;\;\;\;\mbox{w.p.1}$$
immediately translated 
(via the strong approximation of
Theorem~13)
to a strong-law-of-large-numbers
(SLLN) result for the waiting times:
$$\frac{1}{n}\log W_n \to R_1(P,Q,D)
         \;\;\;\;\mbox{w.p.1.}$$

In this section we will prove refinements
to the generalized AEP of Section~2.2, 
and in Section~5 we will revisit the applications
of the previous section and use these refinements
to prove corresponding second order asymptotic 
results.

To get some motivation, 
let us consider for a moment 
the simplest version of the 
classical AEP, for an $\iid$ 
process $\Xp$ with distribution 
$P$ on the finite alphabet $A$. 
The AEP here follows by a simple 
application of the law 
of large numbers,
\be
-\frac{1}{n}\log P^n(X_1^n)
=
\frac{1}{n}\sum_{i=1}^n[-\log P(X_i)]
\to H
\label{eq:oldPS}
\ee
where $H$ is the entropy of $P$.
But (\ref{eq:oldPS}) contains 
more information than that:
It says that $-\log P^n(X_1^n)$ 
is in fact equal to the partial sum 
$S_n=\sum_{i=1}^n Z_i$ of the $\iid$
random variables $Z_i=-\log P(X_i)$.
Therefore we can apply the
central limit theorem (CLT)
or the law of the iterated 
logarithm (LIL) to get more 
precise information on the 
convergence of the AEP.

The same strategy can be carried out
for non-$\iid$ processes: Initially
Ibragimov \cite{ibragimov:62}
and then Philipp and Stout \cite{philipp-stout:book}
showed that even when $\Xp$ is a Markov chain,
or, more generally, a weakly dependent 
random process, the quantities $-\log P^n(X_1^n)$ 
can be approximated by the partial sums of
an associated weakly dependent process.
These results have found a number of
applications in lossless data 
compression and related areas 
\cite{kontoyiannis-jtp}\cite{kontoyiannis-97}.

In this and the following section we will
carry out a similar program in the lossy
case. Throughout this section we will 
adopt the notation and assumptions 
of Section~2.2: Let
$\Xp$ be a stationary ergodic
process with first order marginal
$P_1=P$ on $A$, and let $Q$ be
an arbitrary probability measure
on $\Ahat$. Define $\Dmin$ and $\Dav$,
as before (as in equations (\ref{eq:Dmin}) 
and (\ref{eq:Dav})), and 
assume that $\Dmin<\Dav$ so that
the distortion measure $\rho(X,Y)$ is not
essentially constant in $Y$ with positive
probability. We also impose here the
additional assumption that $\rho$ has
a finite third moment:
\be
D_3\bydef
E_{P\times Q}[\rho^3(X,Y)]<\infty.
\label{eq:third}
\ee

The first result of this section 
refines Theorem~1 by giving a more 
precise asymptotic estimate of the
quantity $-\log Q^n(B(X_1^n,D))$ in
terms of the rate-function $R_1(P,Q,D)$
and the empirical measure $\Phatn$ 
induced by $X_1^n$ on $A^n$
$$\Phatn\bydef\frac{1}{n}\sum_{i=1}^n\delta_{X_i}$$
where $\delta_x$ denotes the measure assigning
unit mass to $x\in A$.

\medskip

{\em Theorem~17: \cite{yang-zhang:99}:}
Let $\Xp$ be a stationary ergodic process
with marginal $P$ on $A$, and let $Q$ be
an arbitrary probability measure on $\Ahat.$
Assume that $D_3=E_{P\times Q}[\rho^3(X,Y)]$ 
is finite. Then for any $D\in(\Dmin,\Dav)$:
\be
-\log Q^n(B(X_1^n,D))= nR_1(\hat{P}_n,Q,D)+\frac{1}{2}\log n + O(1)
\;\;\;\;\mbox{w.p.1.}
\label{eq:br}
\ee

\medskip

Next we show that the most significant
term in (\ref{eq:br}) can be approximated 
by the partial sum of a weakly dependent
random process. Recall the definition of 
the $\alpha$-mixing coefficients of $\Xp$
$$\alpha(k)\;=\;\sup\{|\BBP(A\cap B)-\BBP(A)\BBP(B)|\;:\;\;
A\in\sigma(X_{-\infty}^0),\; B\in\sigma(X_{k}^{\infty})\}$$
where $\sigma(X_i^j)$ is
the $\sigma$-field generated by $X_i^j$. 
The process $\Xp$ is called {\em $\alpha$-mixing}
if $\alpha(k)\to 0$ as $k\to\infty$;
see \cite{bradley} for more details.

We also need to recall some of the notation 
from the proof of Theorem~1 in Section~2.2.
For $x\in A$ and $\la\in\RL$, let $\LA_x(\la)$
denote the log-moment generating function
of the random variable $\rho(x,Y)$
$$\LA_x(\la)\bydef \log E_Q\left(e^{\lambda\rho(x,Y)}\right)$$
and note that the function $\LA(\la)$ defined
in (\ref{eq:GEcheck}) can be written
as $\LA(\la)=E_P[\LA_X(\la)]$.
Also recall that for any $D\in(\Dmin,\Dav)$ there
exists a unique $\la^*<0$ such that 
$\LA'(\la^*)=D$. 


\newpage

{\em Theorem~18: \cite{dembo-kontoyiannis}:}
Let $\Xp$ be a stationary $\alpha$-mixing process
with marginal $P$ on $A$, and let $Q$ be
an arbitrary probability measure on $\Ahat.$
Assume that the $\alpha$-mixing coefficients
of $\Xp$ satisfy
\be
\sum_{k=1}^\infty \alpha^t(k)<\infty, 
\;\;\;\;\mbox{for some $t\in(0,1/3)$}
\label{eq:LIL-cond}
\ee
and that $D_3=E_{P\times Q}[\rho^3(X,Y)]$ 
is finite. Then for any $D\in(\Dmin,\Dav)$:
$$nR_1(\hat{P}_n,Q,D) = nR_1(P,Q,D) + \sum_{i=1}^n 
	g(X_i) + O(\log\log n)
\;\;\;\;\mbox{w.p.1}$$
where 
\be
g(x)\bydef \LA(\la^*) -\LA_x(\la^*),\;\;\;\;
x\in A.
\label{eq:functiong}
\ee

\medskip

Theorem~18 is a small generalization 
of \cite[Theorem~3]{dembo-kontoyiannis}.
Before giving its proof outline,
we combine Theorems~17 and~18 to 
show that, as promised, $-\log Q^n(B(X_1^n,D))$
can be accurately approximated as the 
partial sum of the weakly dependent
random process $\{g(X_n)\}$.

\medskip

{\em Corollary~19: Second Order Generalized AEP:}
Let $\Xp$ be a stationary $\alpha$-mixing process
with marginal $P$ on $A$, and let $Q$ be
an arbitrary probability measure on $\Ahat.$
Assume that the $\alpha$-mixing coefficients
of $\Xp$ satisfy
(\ref{eq:LIL-cond})
and that $D_3=E_{P\times Q}[\rho^3(X,Y)]$ is
finite. Then for any $D\in(\Dmin,\Dav)$, and
with $g(x)$ defined as in (\ref{eq:functiong}):
$$
-\log Q^n(B(X_1^n,D))= nR_1(P,Q,D) + \sum_{i=1}^ng(X_i) 
	+ \frac{1}{2}\log n + O(\log\log n)
\;\;\;\;\mbox{w.p.1.}$$

\medskip

{\em Proof Outline for Theorem~18:}
Adapting the argument leading from
(22) to (24) of \cite{dembo-kontoyiannis},
one easily checks that the result of 
Theorem~18 holds as soon as 
\be
\liminf_{n \to \infty} \inf_{|\theta|<\delta} B_n(\theta)
&>&0
	\;\;\;\;\mbox{w.p.1}
	\label{Bncond}\\
\mbox{and}\;\;\;\;
\limsup_{n \to \infty} \frac{ n A_n^2}{\log \log n}
&<&\infty 
	\;\;\;\;\mbox{w.p.1}
	\label{Ancond}
\ee
where
$A_n=
n^{-1}
\sum_{k=1}^n\zeta_k$ is
the
empirical mean of the centered
random variables $\zeta_k=\Lambda_{X_k}'(\la^*)-D$,
and $B_n(\theta)$ is the
empirical mean
of the non-negative random variables
$\Lambda_{X_k}''(\la^*+\theta)$. 
By the ergodic theorem we have,
with probability one,
\begin{eqnarray*}
\liminf_{n \to \infty} \inf_{|\theta|<\delta} B_n(\theta)
&\geq&
\liminf_{n \to \infty} \frac{1}{n} \sum_{k=1}^n
\inf_{|\theta|<\delta} \Lambda_{X_k}''(\la^*+\theta) \\
&=& E_P \left [ \inf_{|\theta|<\delta} \Lambda_{X}''(\la^*+\theta)
	\right]
\end{eqnarray*}
and by Fatou's lemma and the continuity of
the map $\theta \mapsto \Lambda_{x}''(\la^*+\theta)$
it follows that
$$
\liminf_{\delta \downarrow 0}
E_P \left[
	\inf_{|\theta|<\delta} \Lambda_X''(\la^*+\theta)
	\right]
\geq E_P [\Lambda_X''(\la^*)] = \Lambda''(\la^*) > 0.
$$
This implies that 
(\ref{Bncond}) holds once $\delta>0$
is made small enough. [Note that the above 
argument also avoids an incorrect -- but 
also unnecessary -- application
of the uniform ergodic theorem in the 
derivation of \cite[eq.~(26)]{dembo-kontoyiannis}.]

Turning to (\ref{Ancond}), since $\la^*<0$, 
it follows by the convexity of $\Lambda_x(\la)$ 
that that for any $x\in A$:
$$
0 \leq \Lambda_x'(\la^*) \leq \Lambda_x'(0) = E_Q[\rho(x,Y)].
$$
Consequently, H\"older's inequality and assumption
(\ref{eq:third}) imply that the random variable
$$|\zeta_k| \leq E_Q[\rho(X_k,Y)|X_k]+D$$
has a finite third moment.
Recall  \cite{oodaira-yoshihara:71a}
that the LIL holds for the partial sum $A_n$ of a
zero-mean, stationary process $\{\zeta_k\}$ with
a finite third moment, as soon as 
the $\alpha$-mixing coefficients 
of $\{\zeta_k\}$ satisfy (\ref{eq:LIL-cond}).
The observation 
that $\zeta_k$ is a deterministic
function of $X_k$ for all $k$
completes the proof. \qed

\section{Applications -- Second Order Results}
Here we revisit the applications considered
in Section~3, and using the 
``second order generalized AEP''
of Corollary~19 we prove second order 
refinements for many of the results from 
Section~3. In Section~5.1 we consider
the problem of lossy data compression
in the same setting as in Section~3.1.
We use the second order AEP
to determine the precise asymptotic
behavior of the Shannon random codebooks,
and show that, with probability one,
they achieve optimal compression performance 
up to terms of order $(\log n)$ bits. 
Moreover, essentially the same compression 
performance can be achieved universally. 
For arbitrary variable-length codes 
operating at a fixed rate level, we show 
that the rate at which they can achieve
the optimal rate of $n\calR(D)$ bits is
at best of order $O(\sqrt{n})$ bits. 
This is the 
best possible redundancy rate as 
long as the ``minimal coding variance'' 
of the source is strictly positive. 
For discrete $\iid$ sources, 
a characterization is given of 
when this variance can be zero.

In Section~5.2 we look at waiting times,
and we prove a second order refinement to
Theorem~14, and in Section~5.3 we 
consider the problem of determining
the asymptotic behavior of longest
match-lengths. As discussed briefly
in Section~3.4, their asymptotics
can be deduced from the corresponding
waiting-times results via duality.

\subsection{Lossy Data Compression}

\subsubsection{Random Codes and Second Order Converses}
Here we consider the exact same setup as in Section~3.1:
An $\iid$ source $\Xp$ with distribution $P$ on $A$
is to be compressed with distortion $D$ or less with
respect to a bounded distortion measure
$\rho$, satisfying, as before, the usual
assumption (\ref{eq:maximin}) --
see the remark at the end of this
section for
its implications.
We take 
the reproduction alphabet $\Ahat$ to be
finite, define $\Dbar$ as in (\ref{eq:Dbar}),
and assume that $\Dbar>0$.

For $D\in(0,\Dbar)$, let $Q_n^*$, $n\geq 1$,
denote the optimal reproduction distributions
at distortion level $D$. Combining the
strong approximation Theorem~9 with the
second order generalized AEP
of Corollary~19 and the discussion in 
Section~3.1 yields:


\newpage

{\em Theorem~20: Pointwise Redundancy for I.I.D. Sources
\cite{kontoyiannis-red:00}:}
Suppose $\Xp$ is an $\iid$ source with distribution
$P$ on $A$, and with rate-distortion
function $\calR(D)$ (in bits). Let $Q_n^*$ 
denote the optimal reproduction distributions 
at distortion level $D\in(0,\Dbar)$, 
and define the function
$h(x)=(\log_2 e)g(x)$, $x\in A$,
with $g$ defined as in (\ref{eq:functiong}).
Then:
\begin{itemize}
\item[(a)]
The codes based on almost any realization
of the Shannon random codebooks according
to the measures $\{Q_n^*\}$ have codelengths
$\ell_n(X_1^n)$
satisfying
$$
\ell_n(X_1^n)\leq
n\calR(D)
+\sum_{i=1}^n h(X_i)
+4\log n
\;\;\;\;\mbox{bits, eventually, w.p.1.}$$
\item[(b)]
The codes based on almost any realization 
of the universal Shannon random codebooks 
have codelengths $\ell_n(X_1^n)$ satisfying
$$
\ell_n(X_1^n)\leq
n\calR(D)
+\sum_{i=1}^n h(X_i)
+(4+|\Ahat|)\log n
\;\;\;\;\mbox{bits, eventually, w.p.1.}$$
\end{itemize}

\medskip

We remark that the coefficients of the 
$(\log n)$ terms in (a) and (b) above 
are not the best possible, and can be 
significantly improved; see 
\cite{konto-zhang:00} for more details.  

Perhaps somewhat surprisingly, 
it turns out that the performance 
of the above random codes is 
optimal up to terms of order 
$(\log n)$ bits. 
Recall that a {\em code $C_n$ operating 
at distortion level $D\geq 0$} is
defined by a triplet $(B_n,\phi_n,\psi_n)$ where:
\begin{itemize}
\item[$(a)$]
$B_n$ is a subset of $\Ahatn$, called the {\em codebook},
\item[$(b)$]
$\phi_n:A^n\to B_n$ is the {\em encoder},
\item[$(c)$]
$\psi_n:B_n\to \{0,1\}^*$ is a
uniquely decodable map,
\end{itemize}
such that 
$$\rho_n(x_1^n,\phi_n(x_1^n))\leq D,
\;\;\;\;\;\;\mbox{for all}\;\;x_1^n\in A^n.$$
The codelengths $\ell_n(X_1^n)$ achieved by
such a code are simply:
$$\ell_n(x_1^n)=\;
\mbox{length of}\;[\psi_n(\phi_n(x_1^n))]
\;\;\;\;\mbox{bits}.$$

\medskip

{\em Theorem~21: Pointwise Converse for I.I.D. Sources
\cite{kontoyiannis-red:00}:}
Let $\Xp$ be an $\iid$ source with distribution
$P$ on $A$, and let $\{C_n\}$ be an arbitrary
sequence of codes operating at distortion
level $D\in(0,\Dbar)$, with associated
codelengths $\{\ell_n\}$. Then:
$$
\ell_n(X_1^n)\geq
n\calR(D)
+\sum_{i=1}^n h(X_i)
-\log n
\;\;\;\;\mbox{bits, eventually, w.p.1}$$
where $h(x)$
is defined as in Theorem~20.

\medskip

The proof of Theorem~21 in \cite{kontoyiannis-red:00}
uses techniques quite different to those developed in
this paper. In particular, the key step in the proof 
is established by an application of
the generalized Kuhn-Tucker conditions of Bell and
Cover \cite{bell-cover:88}.

Theorems~20 and~21 are next combined to 
yield ``second order'' refinements to 
Shannon's classical source coding theorem. 
For a source $\Xp$ as in Theorem~21 and 
a $D\in(0,\Dbar)$, the {\em minimal coding
variance $\sigma^2=\sigma^2(P,D)$ of
source $P$ at distortion level $D$}
is
\be
\sigma^2=\sigma^2(P,D)\bydef\VAR[h(X_1)]
\label{eq:mincv}
\ee
with $h(x)$ as in Theorem~20.


\newpage

{\em Theorem~22: Second Order Source Coding Theorems
\cite{kontoyiannis-red:00}:}
Let $\Xp$ be an $\iid$ source with 
distribution $P$ on $A$ and with
rate-distortion function $\calR(D)$
(in bits).
For $D\in(0,\Dbar)$:
\begin{itemize}
\item[]{\bf (CLT)}
There is a sequence of
random variables $G_n=G_n(P,D)$ such that, for any
sequence of codes $\{C_n,\ell_n\}$ operating
at distortion level $D$, we have
\be
\ell_n(X_1^n)-
n\calR(D)\geq \sqrt{n}G_n
        \;\;\;\;\mbox{bits, eventually, w.p.1}
\label{eq:clt}
\ee
and the $G_n$ converge in distribution
to a Gaussian random variable
$$G_n\weakly N(0,\sigma^2)$$
where $\sigma^2=\sigma^2(P,D)$ 
is the minimal coding variance.
\item[]{\bf (LIL)}
With $\sigma^2$ as above,
for any sequence of codes
$\{C_n,\ell_n\}$ operating
at distortion level $D$:
\ben
\limsup_{n\to\infty}\;
\frac{\ell_n(X_1^n)-n\calR(D)}{\sqrt{2n\log\log n}}
&\geq& \sigma\;\;\;\;\mbox{w.p.1}\\
\liminf_{n\to\infty}\;
\frac{\ell_n(X_1^n)-n\calR(D)}{\sqrt{2n\log\log n}}
&\geq& -\sigma\;\;\;\;\mbox{w.p.1.}
\een
\item[]{\bf (\boldmath$\Rightarrow$)}
Moreover, there exist codes $\{C_n,\ell_n\}$
operating at distortion level $D$, that 
asymptotically achieve equality
{\em universally} in all these 
lower bounds.
\end{itemize}

\medskip

{\em Remark on Assumption (\ref{eq:maximin}):}
When the distortion measure does not satisfy 
assumption (\ref{eq:maximin}) [as, for example, 
when $\rho(x,y)=(x-y)^2$ with $A=\RL$ and $\Ahat$
a finite subset of $\RL$], we can modify $\rho$
to $\rho'(x,y)=\rho(x,y)-f(x)$, with 
$f(x)=\min_{y \in \hat{A}} \, \rho(x,y)$, 
so that $\rho'$ satisfies (\ref{eq:maximin}).
Then, to generate codes operating at 
distortion level $D$ with respect to $\rho$, 
we can construct random codebooks for 
as before but do the encoding with respect
to $\rho'(x,y)$ at the {\it random} 
distortion level $D_n= D - E_{\hat{P}_n}(f(X))$. 
It is not hard to check that
\cite[Theorem 2]{dembo-kontoyiannis}
can be extended to apply when $D$ is 
replaced by the sequence $\{D_n\}$.
Since $D_n \to D - E_P(f(X))$ as $n\to\infty$,
this results with the first order 
approximation 
$$-\frac{1}{n}\log Q^*_n(B(X_1^n,D_n))\approx
R_1^{\rho'}(\hat{P}_n,Q^*,D_n).$$
Simple algebra then shows that 
$$
R_1^{\rho'}(\hat{P}_n,Q^*,D_n)=R_1^\rho(\hat{P}_n,Q^*,D)
$$
implying that all the results of Section 5.1.1
remain valid [despite the fact that $\rho$
does not satisfy  (\ref{eq:maximin})], with
the function $h(\cdot)$ taken in terms of 
the log-moment generating function 
$\Lambda_x(\la)$ of the {\it original} 
distortion measure $\rho$ (and not that of 
the modified $\rho'$).

\subsubsection{Critical Behavior}
In view of Theorems~20 and~21 above,
the codelengths $\ell_n^*(X_1^n)$ of 
the best code operating at distortion 
level $D$ have:
$$\ell^*_n(X_1^n)\approx
n\calR(D)
+\sum_{i=1}^n h(X_i)
+O(\log n)
\;\;\;\;\mbox{bits.}$$
This reveals an interesting 
dichotomy in the behavior of 
the ``pointwise'' redundancy of
the best code:
\begin{itemize}
\item
Either the minimal coding variance
$\sigma^2$ (recall (\ref{eq:mincv})) 
is nonzero, in which case the best
rate at which optimality can
be achieved is of order $\sqrt{n}$
bits by the CLT;
\item
or $\sigma^2=0$, and the best redundancy 
rate is of order $(\log n)$ bits
(cf. \cite{zhang-yang-wei:I}).
\end{itemize}
Under certain conditions, in this section 
we give a precise characterization of 
when each of these two cases can occur. 
Before stating it, we briefly discuss 
two examples to gain some intuition.
 
\medskip

{\em Example~4: Lossless Compression:}
Lossless data compression 
can be considered as an extreme case 
of lossy compression, where $\Xp$ is
an $\iid$ source with distribution $P$
on a finite set $A=\Ahat$,
and the distortion level $D$ 
is set to zero. Here it is 
well-known that (ignoring the integer 
length constraints) the best code is 
given by the idealized Shannon code,
$\ell_n(X_1^n)=-\log_2 P^n(X_1^n)$.
Accordingly, the upper 
and lower bounds of Theorems~21 
and~22 say that the best code has 
codelengths 
$$\ell_n(X_1^n) = n\calH(P)
	+\sum_{i=1}^n h(X_i)$$
where $\calH(P)$ is the entropy of $P$
in bits, and with
$$h(x)\bydef-\log_2 P(x) - \calH(P),
	\;\;\;\;x\in A.$$
When is $\sigma^2=0$? By its
definition (\ref{eq:mincv}),
$\sigma^2$ is zero if and only if
the function $h(x)$ is constant over $x$,
which, in this case, can only happen if
$P(x)$ is constant over $x\in A$. 
Therefore, here:
{\em $\sigma^2=0$ if and only if
the source has a uniform distribution
over $A$.} 

\medskip
 
{\em Example~5: Binary Source with Hamming Distortion:}
Consider the simplest
non-trivial lossy example: 
Let $\Xp$ be an $\iid$ source
with Bernoulli($p$) distribution 
(for some $p\in(0,1/2]$),
let $A=\Ahat=\{0,1\}$,
and take $\rho$ to be Hamming
distortion: $\rho(x,y)=|x-y|$.
For $D\in(0,p)$ it is not
hard to evaluate all the 
relevant quantities 
explicitly
(see, e.g., 
\cite[Example~2.7.1]{berger:book}
or \cite[Theorem~13.3.1]{cover:book}).
In particular,
the optimal reproduction
distribution $Q^*$ is
Bernoulli($q$),
with $q=(p-D)/(1-2D)$, and
our function of interest is:
$$h(x)=
-\log_2\left(\frac{P(x)}{1-D}\right)
        -E_P\left[
                -\log_2\left(\frac{P(X_1)}{1-D}\right)
                \right].$$
Recalling that the minimal coding
variance is zero if and only if
$h(x)$ is constant, from the above
expression we see that, similarly
to the previous example, also 
here:
{\em $\sigma^2=0$ if and only if
the source has a uniform distribution}.

\medskip

For discrete sources, the next result gives
conditions under which the characterization 
suggested by these two examples remains valid.
Suppose $A=\Ahat=\{a_1,a_2,\ldots,a_k\}$
is a finite set, write $\rho_{ij}$ for 
$\rho(a_i,a_j)$, and assume
that $\rho$ is symmetric
and that $\rho_{ij}=0$ if and only if
$i=j$. We call $\rho$ a 
{\em permutation distortion measure}, 
if all rows of the matrix 
$(\rho_{ij})_{i,j=1,\ldots,k}$
are permutations of one another.

\medskip

{\em Theorem~23: Variance Characterization
\cite{dembo-kontoyiannis:crit:01}:}
Let $\Xp$ be a discrete source with
distribution $P$ and rate-distortion 
function $R(D)$. Assume that $R(D)$ 
is strictly convex over $(0,\Dbar)$. 
There are exactly two possibilities:
\begin{itemize}
\item[(a)]
Either $\sigma^2=\sigma^2(P,D)$ is only
zero for finitely many $D\in(0,\Dbar).$
\item[(b)]
Or $\sigma^2=\sigma^2(P,D)\equiv 0$ 
for {\em all} $D\in(0,\Dbar)$, in which 
case $P$ is the uniform distribution
on $A$ and $\rho$ is
a permutation distortion measure.
\end{itemize}

\medskip

A general discussion of this
problem, including the case of continuous
sources, is given in 
\cite{dembo-kontoyiannis:crit:01}.
Also, in the lossless case,
the problem of characterizing 
when $\sigma^2=0$ for sources 
with memory is dealt with 
in \cite{kontoyiannis-97}.

Before moving on to waiting times and match-lengths 
we mention that, in a somewhat similar
vain, the problem of understanding the best
{\em expected}
redundancy rate in lossy data
compression has also been recently considered in 
\cite{zhang-yang-wei:I})\cite{yang-zhang:II}%
\cite{yang-zhang:III}\cite{ishii-yamamoto:97}.

\subsection{Waiting Times}
Next we turn to waiting times.
Recall that, given $D\geq 0$ 
and two independent realizations 
of the stationary ergodic 
processes $\Xp$ and $\Yp$,
the waiting time $W_n$ was
defined as the time of the
first appearance of $X_1^n$
in $\Yp$ with distortion $D$
or less (see (\ref{eq:Wn-def})
for the precise definition).
In Theorem~14 we gave conditions
that identified the first order
limiting behavior of $W_n$.
In particular, when $\Yp$ is
$\iid$, it was shown in 
Theorem~14~(a)
that
\be
\frac{\log W_n}{n}\to R_1(P,Q,D)
\;\;\;\;\mbox{w.p.1}
\label{eq:w-slln3}
\ee
where $P$ and $Q$ are the first
order marginals of $\Xp$ 
and $\Yp$, respectively.

The next result gives conditions
under which the SLLN-type
statement of (\ref{eq:w-slln3})
can be refined to a CLT and
a LIL.

\medskip

{\em Theorem~24: CLT and LIL for Waiting Times:}
Let $\Xp$ be a stationary $\alpha$-mixing process
and $\Yp$ be an $\iid$ process, with marginal
distributions $P$ and $Q$, on $A$ and $\Ahat$,
respectively. Assume that the $\alpha$-mixing 
coefficients of $\Xp$ satisfy (\ref{eq:LIL-cond})
and that $D_3=E_{P\times Q}[\rho^3(X,Y)]$ is
finite. Then for any $D\in(\Dmin,\Dav)$ the
following series converges
\be
\sigma^2\bydef E_P[g^2(X_1)]+2\sum_{k=2}^\infty E_P[g(X_1)g(X_k)]
\label{eq:variance}
\ee
with $g(x)$ defined as in (\ref{eq:functiong}),
and, moreover:
\begin{itemize}
\item[]{\bf (CLT)} With $R_1=R_1(P,Q,D)$:
$$\frac{\log W_n \;-\; nR_1}{\sqrt{n}}
	\weakly N(0,\sigma^2).$$
\item[]{\bf (LIL)}
The set of limit points of the sequence
$$\left\{
	\frac{\log W_n \;-\; nR_1}
	     {\sqrt{2n\log\log n}}
  \right\},\quad n\geq 3$$
coincides with $[-\sigma,\sigma]$, with
probability one.
\end{itemize}

\medskip

{\em Proof Outline:}
For a bounded distortion measure 
$\rho$, Theorem~24 was proved in 
\cite{dembo-kontoyiannis}. 
To obtain the more general statement above
combine the strong approximation
of Theorem~13 with the second order
AEP in Corollary~19 to get:
\be
\log W_n=
nR_1(P,Q,D) + \sum_{i=1}^ng(X_i) + O(\log n)
\;\;\;\;\mbox{w.p.1.}
\label{eq:inter}
\ee
Since $\Xp$ satisfies the mixing
assumption (\ref{eq:LIL-cond}),
so does the process $\{g(X_n)\}$.
Also, since $\la^*<0$, the function
$\LA_x(\la^*)$ is bounded above by zero,
and by Jensen's inequality it is
bounded below by $\la^*E_Q[\rho(x,Y)].$
Therefore,
$$|\LA_x(\la^*)|\leq |\la^*|E_Q[\rho(x,Y)]$$
and this, together with 
H\"older's inequality and
the definition of $g(x),$ imply
that $E_P[|g(X_1)|^3]<\infty$.
Therefore we can apply the CLT
of \cite[Theorem~1.7]{peligrad:86} 
to the process $\{g(X_n)\}$
in order to deduce the CLT-part 
of the theorem from (\ref{eq:inter}).
Similarly, applying the LIL of
\cite{oodaira-yoshihara:71a}
to $\{g(X_n)\}$, from (\ref{eq:inter})
we get the LIL-part of the theorem.
\qed

\medskip

{\em Remark 5:} When the variance
$\sigma^2$
in (\ref{eq:variance}) is positive,
then the {\em functional} versions of 
the above CLT and LIL given in 
\cite{dembo-kontoyiannis} still hold, 
under exactly the conditions of Theorem~24.
(This follows by
applying the functional CLT of
\cite[Theorem~1.7]{peligrad:86}
and the functional LIL of
\cite[Theorem~1~(IV)]{oodaira-yoshihara:71b}.)

\subsection{Match-Lengths and Duality}

Finally we turn to our last application,
match-lengths. Recall that, given a 
distortion level $D\geq 0$ and two 
independent realizations of the 
processes $\Xp$ and $\Yp$, the match-length 
$L_m$ is defined as the length $\ell$ 
of the longest prefix $X_1^\ell$ that 
appears (with distortion $D$ or less) 
starting somewhere in the ``database'' 
$Y_1^m.$ See (\ref{eq:Lm-def}) for the 
precise definition. As we briefly mentioned
in Section~3.4, there is a duality
relationship between match-lengths
and waiting times: Roughly speaking,
long matches mean short waiting times,
and vice-versa;
see (\ref{eq:duality}).

Although the relation (\ref{eq:duality}) 
is not as simple as the duality 
(\ref{eq:easy-dual}) for exact matching,
it is still possible to use 
(\ref{eq:duality}) to translate
the asymptotic results for $W_n$
to corresponding results for $L_m$.
These are given in Theorem~25 below.
This translation, carried out
in \cite{dembo-kontoyiannis}, is 
more delicate than in the case of
exact matching. For example, in
order to prove the CLT for the 
match-lengths $L_m$ one
invokes
the functional CLT for
the waiting times (see Remark~5 above
and the proof of Theorem~4 in 
\cite{dembo-kontoyiannis}).

\medskip

{\em Theorem~25: Match-Lengths Asymptotics:}
Let $\Xp$ be a stationary process
and $\Yp$ be an $\iid$ process, with marginal
distributions $P$ and $Q$, on $A$ and $\Ahat$,
respectively. Assume 
that $D_3=E_{P\times Q}[\rho^3(X,Y)]$ is
finite. Then for any $D\in(\Dmin,\Dav)$ we have
$$\mbox{\bf (LLN)}\hspace{1.65in} 
\frac{L_m}{\log m}\,\to\,\frac{1}{R_1}\;\;\;\;
\mbox{w.p.1}\hspace{2.4in}$$
where $R_1=R_1(P,Q,D)$.
If, moreover, 
the $\alpha$-mixing coefficients of $\Xp$
satisfy (\ref{eq:LIL-cond}) and 
the variance $\sigma^2$ in (\ref{eq:variance}) 
is nonzero, then, with $\tau^2\bydef \sigma^2R_1^{-3}$,
we have,
\ben
&\mbox{\bf (CLT)}&
	\hspace{1.6in}
	\frac{L_m-\frac{\log m}{R_1}}{\sqrt{\log m}}
	\,\weakly\,N(0,\tau^2)
	\hspace{2.2in}\\
&\mbox{\bf (LIL)}&
	\hspace{1.2in} 
	\limsup_{
	m\to\infty}
	\,\frac{L_m-\frac{\log m}{
R_1
}}
	{\sqrt{2\log m\,\log\log\log m}}\,
	=\,\tau\;\;\;\;
	\mbox{w.p.1.}
\een

\medskip

The results of Theorem~25 were
proved in \cite{dembo-kontoyiannis}
for any bounded distortion measure 
$\rho$.  The slightly 
more general version stated above
is proved in exactly the same way, 
using the results of Section~4 
in place of Theorems~2 and~3 
of \cite{dembo-kontoyiannis}.

\section{Random Fields -- First Order Results}
This and the following section are devoted 
to generalizations of the results of
Sections~2--5 to the case of random fields.
Specifically, the role of the processes $\Xp$
and $\Yp$ will now be played by stationary
ergodic random fields
$\Xp=\{X_u\;;\;u\in\IN^d\}$
and $\Yp=\{Y_u\;;\;u\in\IN^d\}$.
As we will see, many of the problems 
that we considered have natural 
analogs in this case, and 
the overall theme
carries over:
The generalized AEP and its refinement
can be extended to random fields,
and the corresponding questions in
data compression and pattern matching
can be answered following 
the same path as before.

\subsection{Notation and Definitions}
The following definitions and notation
will remain in effect throughout Sections~6
and~7.

We consider two random fields
$\Xp=\{X_u\;;\;u\in\IN^d\}$
and $\Yp=\{Y_u\;;\;u\in\IN^d\}$,
$d\geq 2$, taking values 
in $A$ and $\Ahat$, 
respectively, and indexed
by points $u=(u_1,u_2,\ldots,u_d)$
on the integer lattice $\IN^d$.
As before, $A$ and $\Ahat$
are complete, separable 
metric spaces, equipped with 
their Borel
$\sigma$-fields ${\cal A}$
and $\hat{\cal A}$, 
respectively.
Let $\BBP$ and $\BBQ$
denote the 
(infinite-dimensional)
measures of the entire random
fields $\Xp$ and $\Yp$.
Unless explicitly stated
otherwise, we always assume
that $\Xp$ and $\Yp$ are
independent of each other.

Throughout the rest of the
paper we will
assume that $\Xp$ and $\Yp$
are stationary and ergodic. 
To be precise, by that we mean 
that the Abelian group of 
translations
$\{T_u\,:\,u\in\IN^d\}$
acts on both 
$(A^{\IN^d},{\cal A}^{\IN^d},\BBP)$
and
$(\hat{A}^{\IN^d},\hat{\cal A}^{\IN^d},\BBQ)$
in a measure-preserving,
ergodic manner; see \cite{krengel:book}
for a detailed exposition.

For $v,w\in\IN^d$,
the distance between $v$ and $w$ 
is defined by
$$d(v,w)\bydef\max_{1\leq i\leq d}|v_i-w_i|$$
and the distance between two subsets 
$V,W\subset\IN^d$ is
$$d(V,W)\bydef\inf_{v\in V,\;w\in W} d(v,w).$$
Given $v,w\in\IN^d$, we let
$[v,w]=\{u\in\IN^d\;:\;
\mbox{$v_j\le u_j\leq w_j$ for all $j$}\}$,
where $[v,w]$ is empty in case $v_j>w_j$ for some $j$.

We write $C(n)$ for the 
$d$-dimensional cube of side $n\geq1$,
\ben
C(n)=
	\{u\in\IN^d\;:\;\mbox{$1\leq u_j\leq n$ for all $j$}\}
\een
and $[0,\infty)$ for the ``infinite cube''
$$[0,\infty)=\{u\in\IN^d\;:\;
	\mbox{$u_j\geq 0$ for all $j$}\}.$$
For an arbitrary subset 
$U\subset\IN^d$ we let
$|U|$ denote its size;
for example, $|C(n)|=n^d$.
Also for $U\subset\IN^d$ we write
$$X_U\bydef\{X_u\;;\;u\in U\}$$
so that, in particular,
$X_{[0,\infty)} = \{X_u\;;\;
	\mbox{$u_j\geq 0$ for all $j$}\}.$
For 
$V\subset\IN^d$ and $u\in\IN^d$ we
let $u+U$ denote the translate
$$u+V=\{u+v\;:\;v\in V\}.$$

For each $n\geq 1$, let $P_n$ denote the 
marginal distribution of $X_{C(n)}$
on $A^{n^d}$, and similarly write
$Q_n$ for the distribution of $Y_{C(n)}$.
Let $\rho:A\times\Ahat\to[0,\infty)$
be an arbitrary nonnegative (measurable)
function, and define a sequence of
single-letter distortion measures 
$\rho_n:A^{n^d}\times\Ahatnd\to[0,\infty)$,
$n\geq 1$, by
\ben
\rho_n(x_{C(n)},y_{C(n)})\bydef\frac{1}{n^d}
	\sum_{u\in C(n)}\rho(x_u,y_u)
\;\;\;\;x_{C(n)}\in A^{n^d},\;y_{C(n)}\in\Ahatnd.
\een
Given $D\geq 0$ and $x_{C(n)}\in A^{n^d}$,
we write 
$B(x_{C(n)},D)$ for
the distortion-ball of radius $D$:
$$B(x_{C(n)},D)=
\left\{
y_{C(n)}\in\Ahatnd\;:\;\rho_n(x_{C(n)},y_{C(n)})\leq D
\right\}.$$

\subsection{Generalized AEP}
It is well-known that the classical AEP
\ben
-\frac{1}{n}\log P_n(X_1^n) \to H(\BBP)
\;\;\;\;\mbox{w.p.1}
\een
generalizes to the case of finite-alphabet 
random fields on $\IN^d$, as well 
as to other amenable group actions
\cite{ornstein-weiss:83}. In this 
section we give two versions of the
generalized AEP of Theorems~1 and~4
to the case of random fields on $\IN^d$.

\paragraph{$\Yp$ is i.i.d. }
In the notation of Section~6.1,
we take $\Xp$ to be a stationary 
ergodic random field with first 
order marginal $P_1=P,$ and 
$\Yp$ to be i.i.d. with first 
order marginal $Q_1=Q$. 
We define $\Dmin$ and $\Dav$ 
as in the one-dimensional 
case (recall equations (\ref{eq:Dmin})
and (\ref{eq:Dav})), and assume
that $\rho(x,y)$ is not essentially 
constant for ($\BBP$-almost) 
all $x\in A$, that is, $\Dmin < \Dav.$

A simple examination of the proof of 
Theorem~1 shows that it
extends {\sl verbatim} to the
case of random fields, with the
only difference that instead of the
usual ergodic theorem we now need
to invoke the ergodic theorem
for $\IN^d$ actions; see
\cite[Chapter~6]{krengel:book}.
We thus obtain:

\medskip

{\em Theorem~26. Generalized AEP when $\Yp$ is $\iid$:}
Let $\Xp$ be a stationary ergodic random field on
$\IN^d$ and $\Yp$ be $\iid$, with marginal distributions
$P$ and $Q$ on $A$ and $\Ahat$, respectively.
Assume that $\Dav=E_{P\times Q}[\rho(X,Y)]$ is
finite. Then for any $D\in(\Dmin,\Dav)$
\ben
-\frac{1}{n^d}\log Q^{n^d}(B(X_{C(n)},D)) \to R_1(P,Q,D)
        \;\;\;\;\mbox{w.p.1}
\een
with the (one-dimensional)
rate-function $R_1(P,Q,D)$
defined as in Theorem~1.


\paragraph{$\Yp$ is not i.i.d. }
Let $\Xp$ and $\Yp$ be stationary random fields and
define $\Dav$ and $\Dmax$ exactly as in the
one-dimensional case (recall 
(\ref{eq:Dav}) and (\ref{eq:Dmax})).
We assume that the distortion 
measure $\rho$ is essentially
bounded, $\Dmax < \infty$,
and define
\be
\Dmin \bydef 
\sup_{n \geq 1} \Dminn = \lim_{n\to\infty} \Dminn 
\label{eq:dmind}
\ee
where 
\be
\label{eq:dminn}
\Dminn \bydef 
E_{P_n}[\essinf_{Y_{C(n)}\sim Q_n} \;\rho_n(X_{C(n)},Y_{C(n)})].
\ee
To see that the limit in (\ref{eq:dmind}) 
exists and equals the supremum, first 
note that $\{n^d \Dminn\}$ is an
increasing sequence, and that 
$D_{\rm min}^{(nk)}\geq D_{\rm min}^{(k)}$
for all $n,k\geq 1$.
Now fix $k\geq 1$ arbitrary. Given $n\geq k$
we write $n = mk + r$ for some $0\leq r\leq k-1$,
so that
$$ n^dD_{\rm min}^{(n)}\geq (mk)^dD_{\rm min}^{(mk)}
\geq (mk)^dD_{\rm min}^{(k)}.$$
Since $n/mk\to 1$ as $n\to\infty$, this implies that
$$\liminf_{n\to\infty} D_{\rm min}^{(n)} \geq D_{\rm min}^{(k)}.$$
Since $k$ was arbitrary we are done.

Finally, we assume once again that
the distortion measure $\rho$ is
not essentially constant, 
that is, $\Dmin<\Dav$.
Our next result is the 
random fields analog of Theorem 4; 
it is proved in Appendix~C.

\medskip

{\em Theorem~27. Generalized AEP rate function.}
Let $\Xp$ and $\Yp$ be stationary random fields.
Assume that $\rho$ is bounded, and that with
$\BBP$-probability one, conditional on 
$X_{[0,\infty)}=
x_{[0,\infty)}$, the random variables
$\{\rho_n(x_{C(n)},Y_{C(n)} )\}$ satisfy a
large deviations principle with some
deterministic, convex rate-function.
Then for all $D\in (\Dmin,\Dav)$, 
except possibly at $D = \Dinf$, 
\be
\label{eq:ldp-27}
\lim_{n \to \infty}
-\frac{1}{n^d}\log Q_n(B(X_{C(n)},D)) = R(\BBP,\BBQ,D)
        \;\;\;\;\mbox{w.p.1}
\ee
where $\Dinf$ 
and the rate-function $R(\BBP,\BBQ,D)$ 
are defined as in the one-dimensional case, 
by (\ref{eq:dinf}) and (\ref{eq:thm4b}),
respectively, and the rate-functions
$R_n(P_n,Q_n,D)$ are now defined as
\ben
R_n(P_n,Q_n,D) = \inf_{V_n} \frac{1}{n^d} H(V_n\|P_n\times Q_n)
\een
with the infimum taken over all joint distributions 
$V_n$ on $A^{n^d}\times\Ahatnd$ such that
the $A^{n^d}$-marginal of $V_n$ is $P_n$
and $E_{V_n}[\rho_n(X_{C(n)},Y_{C(n)})]\leq D$.

\medskip

{\em Remark 6:} Suppose that $(\Xp,\Yp)$
is a stationary random field satisfying 
a ``process-level LDP'' with a convex, good
rate-function. To be precise,
given $x_{C(n)}\in A^{n^d}$, 
write $x^{(n)}$ for the periodic
extension of $x_{C(n)}$ to an 
infinite realization in $A^{[0,\infty)}$
and let $X^{(n)}$ and $Y^{(n)}$ denote
the periodic extensions of $X_{C(n)}$
and $Y_{C(n)}$, respectively.
The process-level empirical 
measure $\calLn$  induced
by $\Xp$ and $\Yp$ on 
$(A^{[0,\infty)}\times\hat{A}^{[0,\infty)})$ is
defined by
$$\calLn\bydef\frac{1}{n^d}\sum_{u\in C(n)}
	\delta_{(X^{(n)}_{u+[0,\infty)},Y^{(n)}_{u+[0,\infty)})}$$
where $\delta_{s,s'}$ denotes the measure
assigning unit mass to the joint realization
$(s,s')\in A^{[0,\infty)}\times\hat{A}^{[0,\infty)}$,
and $X^{(n)}_{u+[0,\infty)}$ 
(or $Y^{(n)}_{u+[0,\infty)}$) 
denotes $X^{(n)}$ 
(respectively, $Y^{(n)}$) 
shifted by $u$ [i.e., the
value of $X^{(n)}_{u+[0,\infty)}$ 
at position $v$ is the same as
the value of $X^{(n)}$ at position
$u+v$; similarly for $Y^{(n)}_{u+[0,\infty)}$.]
By assuming that $(\Xp,\Yp)$
satisfy a ``process-level LDP''
we mean that the sequence of measures
$\{\calLn\}$ satisfies the LDP in 
the space of stationary
probability measures on
$(A^{[0,\infty)}\times\hat{A}^{[0,\infty)})$
equipped with the topology of weak convergence,
with some convex, good rate-function $I(\cdot)$.
These assumptions are satisfied by many of 
the random field models used in applications, 
and in particular by a large class of Gibbs fields 
(see, e.g., 
\cite{comets:86}\cite{folmer-orey}\cite{olla:88} 
for general theory and 
\cite{guyon:book}\cite{winkler:book} for
examples in the areas of
image processing and image analysis).

As in the one-dimensional case, suppose
that the process-level LDP condition 
holds, and that the 
distortion measure $\rho$ is
bounded and continuous on $A\times\Ahat$.
Then with $\BBP$-probability one,
conditional on 
$X_{[0,\infty)}=x_{[0,\infty)}$,
the sequence 
$\{\rho_n(x_{C(n)},Y_{C(n)})\}$ satisfies 
the LDP upper bound with respect to the 
deterministic, convex rate-function $J(\cdot)$
as in Remark~3.
Moreover, assuming sufficiently strong mixing 
properties for $\Yp$ one may also verify the 
corresponding lower bound (for example, by 
adapting the stochastic subadditivity approach of 
\cite{chi-AP:01}).


\subsection{Applications}
In Sections~6.3.1 and~6.3.2 below we
consider the random field analogs of 
the problems discussed in Section~3
in the context of one-dimensional 
processes. 
In the instances when our analysis 
was restricted to $\iid$ processes, 
the extension to random fields is 
trivial -- an $\iid$ random field 
is no different from an $\iid$ process. 
For that reason, we only give the 
full statements of corresponding 
random fields results when the 
generalization from $d=1$ to 
$d\geq 2$ does involve some
modifications. Otherwise, only 
a brief description of the corresponding 
results is mentioned.

\subsubsection{Lossy Data Compression}
Here we very briefly discuss the 
problem of data compression, when 
the data is in the form of a two- 
or more generally a $d$-dimensional 
array.
In this case, the underlying data 
source is naturally modeled as 
a $d$-dimensional random field.
Extensive discussions of the
general information-theoretic
problems on random fields are
given 
in \cite{berger-shen-ye:92}
and 
the recent monograph \cite{ye-berger:book};
see also 
\cite{follmer:73}.

First we discuss the results
given in Section~3.1. The construction
of the random codebooks described there
generalizes to random fields in an 
obvious fashion, and the statement 
as well as the proof of Theorem~9
remain unchanged. Following the
notation exactly as developed for 
$\iid$ sources, the strengthened 
coding theorems given in 
Theorems~10 and~11 follow by 
combining (the obvious 
generalization of) Theorem~9 
with the generalized
AEP of Theorem~26.

Similarly, the mismatched-codebook
results of Section~3.2 only rely on
Theorem~9 and the generalized AEP
of Theorem~1, and therefore 
immediately generalize to the
random field case. Finally 
Theorems~15 and~16 in Section~3.5
are only stated for $\iid$ processes,
hence, as mentioned above, they
trivially extend to random 
fields.

\subsubsection{Waiting Times}
Here we consider the natural
$d$-dimensional analogs of the
waiting times questions considered
in Section~3.3. Given two 
independent realizations
of the random fields $\Xp$ and $\Yp$,
our main quantity of interest here
is how ``far'' we have to look
in $\Yp$ until we find a match for
the pattern $X_{C(n)}$ with distortion
$D$ or less. Given $n\geq 1$ and 
a distortion level $D\geq 0$, we 
define the {\em waiting time} $W_n$
as the smallest length $i$ 
such that a copy of the pattern 
$X_{C(n)}$ appears somewhere in 
$Y_{C(i+n-1)}$, with distortion
$D$ or less.
Formally,
\ben
W_n\;=\;\inf\{i\geq 1\; :\; 
	\rho_n(X_{C(n)},Y_{u+C(n)})\leq D
	\;\;\mbox{for some}\;u\in[0,i-1]^d\}
\een
with the convention that the infimum 
of the empty set equals $+\infty$.

In the one-dimensional case our
main tool in investigating
the asymptotic behavior of the
waiting times was the strong
approximation in Theorem~13.
Roughly speaking, Theorem~13
stated that the waiting time
$W_n$ for a $D$-close match
of $X_1^n$ in $\Yp$ is 
inversely proportional
to the probability $Q_n(B(X_1^n,D))$
of such a match.
In Theorem~28 below we generalize
this result to the $d$-dimensional
case by showing that the $d$-dimensional
volume $(W_n)^d$ we have to search 
in $\Yp$ in order to find a $D$-close
match for $X_{C(n)}$ is, roughly,
inversely proportional
to the probability $Q_n(B(X_{C(n)},D))$
of finding such a match.

Before stating
Theorem~28 we need to recall
the following definition.
Dobrushin's {\em non-uniform 
$\phi$-mixing coefficients}
of a stationary random field
$\Yp$ are
\ben
\phi_\ell(k)\;=\;\sup\{|\BBQ(B|A)-\BBQ(B)|\;:
& & \hspace{-0.2in}
	B\in\sigma(Y_{U}),\; A\in\sigma(Y_{V}),\; \BBQ(A)>0\\
& & \quad\quad |U|\leq \ell,\; |V|<\infty,\; d(U,V)\geq k \}
\een
where $\sigma(Y_U)$ denotes
the $\sigma$-field generated by 
the random variables $Y_U$, $U\subset\IN^d$. 
See \cite[Chapter~6]{lin-lu:book}
or \cite{doukhan:book} for detailed
discussions of the coefficients
$\{\phi_\ell(k)\}$ and their properties.

\medskip

{\em Theorem~28. Strong Approximation:}
Let $\Xp$ and $\Yp$ be stationary ergodic 
random fields, and assume that the non-uniform
$\phi$-mixing  coefficients of $\Yp$ satisfy
\be
\limsup_{n\to\infty}\sum_{j=1}^\infty
(j+1)^{d-1}\phi_n(jn)<\infty.
\label{eq:dobrushin}
\ee
If $Q_n(B(X_{C(n)},D))>0$ eventually with 
probability one,
then for any $\epsilon>0$:
\ben
-(1+\epsilon)\log n
\;\leq\;
\log [W^d_n Q_n(B(X_{C(n)},D))]
\;\leq\;
(d+1+\epsilon)\log n
\;\;\;\;\mbox{eventually, w.p.1.}
\een

\medskip

The proof of Theorem~28 is a 
straightforward modification of
the corresponding one-dimensional argument in
\cite{dembo-kontoyiannis}; it is given in 
Appendix~D.

\medskip

{\em Remark 7:} The mixing condition
(\ref{eq:dobrushin}) is satisfied by
a rather large class of 
stationary random fields. For
example in the case of Markov 
random fields, it is easy to check 
that under Dobrushin's uniqueness 
condition the limit in 
(\ref{eq:dobrushin}) is finite;
see \cite[Section~8.2]{georgii:1}
or \cite{doukhan:book} for 
more details.

\medskip

Next we combine the above strong
approximation result with the
generalized AEPs of Theorems~26
and~27, to read off the first order
asymptotic behavior of the 
waiting times. Theorem~29
below generalizes Theorem~14
to the random field case.

\medskip
 
{\em Theorem~29. SLLN for Waiting Times:}
Let $\Xp$ and $\Yp$ be stationary ergodic 
random fields:
 
(a)~If $\Yp$ is $\iid$ and the
average distortion $\Dav$ is finite,
then for any $D\in(\Dmin,\Dav)$
\ben
\frac{1}{n^d}\log W_n^d \to R_1(P_1,Q_1,D)
\;\;\;\;\mbox{w.p.1.}
\een
 
(b)~Suppose that the 
conditions of Theorem~27 are satisfied,
and that $\Yp$ also satisfies 
the mixing assumption (\ref{eq:dobrushin}).
Then, for any $D\in(\Dinf,\Dav)$:
\ben
\frac{1}{n^d}\log W^d_n \to R(\BBP,\BBQ,D)
\;\;\;\;\mbox{w.p.1.}
\een
 
\medskip


\newpage

\section{Random-Fields -- Second Order Results}
Finally we turn to the random field extensions 
of the second order results of Sections~4 and~5.
In Section~7.1 we state the random field analog 
of the second order generalized AEP, and in~7.2
we discuss its application to the problems
of lossy data compression and pattern matching.

\subsection{Refinements of Generalized AEP}
Let $\Xp$ be a stationary ergodic random 
field with marginal distribution $P$ on $A$,
and let $Q$ be a fixed probability measure
on $\Ahat$. We will assume throughout that
the distortion measure $\rho$ has a finite
third moment,
\be
D_3\bydef
E_{P\times Q}[\rho^3(X,Y)]<\infty
\label{eq:third-d}
\ee
and that it is not essentially
constant, i.e., $\Dmin<\Dav$,
with $\Dmin$ and $\Dav$ defined
as before (cf. (\ref{eq:Dmin})
and (\ref{eq:Dav})).

The goal of this section is
to give the random field analogs of 
Theorems~17 and~18 and of Corollary~19
from the one-dimensional case.

An examination of the proof of Theorem~17 in
\cite{yang-zhang:99} shows that its proof
only depends on the ergodicity of $\Xp$
and the $\iid$ structure of the product 
measures $Q^n$. Simply replacing the
application of the ergodic theorem
by the ergodic theorem
for $\IN^d$ actions
\cite[Chapter~6]{krengel:book}
immediately yields the following 
generalization: As long as condition
(\ref{eq:third-d}) is satisfied, 
for all $D\in(\Dmin,\Dav)$ we have
\be
-\log Q^{n^d}(B(X_{C(n)},D))= n^dR_1(\hat{P}_n,Q,D)+\frac{d}{2}\log n + O(1)
\;\;\;\;\mbox{w.p.1}
\label{eq:br-d}
\ee
where $\hat{P}_n$ is now the empirical measure
induced by $X_{C(n)}$ on $A$.

In order to generalize Theorem~18 to $\IN^d$
we need to introduce a measure of dependence
analogous to $\alpha$-mixing in the 
one-dimensional case. For a stationary
random field $\Xp$ on $\IN^d$ we define
the {\em uniform $\alpha$-mixing coefficients}
of $\Xp$ by
\ben
\alpha(k)\;=\;\sup\{|\BBP(A\cap B)-\BBP(A)\BBP(B)|\;:
& & \hspace{-0.2in}
        A\in\sigma(X_{U}),\; B\in\sigma(X_{V}),\; d(U,V)\geq k \}
\een
where, as before, $\sigma(X_U)$ denotes
the $\sigma$-field generated by
the random variables $Y_U$.
See \cite{lin-lu:book}\cite{doukhan:book} 
for more details.

Apart from ergodicity, the main technical
ingredient in the proof of Theorem~18 above
(see also the proof of 
\cite[Theorem~3]{dembo-kontoyiannis})
is the LIL for $\Xp$. 
Similarly to the one-dimensional case, 
the LIL for a random field $\Xp$
holds as soon as the following
mixing condition is satisfied
\be
\alpha(k)\leq C\ k^{
-
3d(1+\epsilon)},
\quad\mbox{for some $\epsilon>0$ and $C<\infty.$}
\label{eq:LIL-cond2}
\ee
[This follows from the 
almost sure invariance principle
in \cite[Theorem~1]{berkes-morrow}.]

Assuming that (\ref{eq:LIL-cond2})
and the third moment condition
(\ref{eq:third-d}) both hold,
we get the following generalization
of Theorem~18. For all $D\in(\Dmin,\Dav)$,
\be
n^dR_1(\hat{P}_n,Q,D) = n^dR_1(P,Q,D) + \sum_{u\in C(n)}
        g(X_u) + O(\log\log n)
\;\;\;\;\mbox{w.p.1}
\label{eq:taylor-d}
\ee
with $g(x)$ defined exactly as in the
one-dimensional case (\ref{eq:functiong}).

Combining (\ref{eq:br-d}) and (\ref{eq:taylor-d})
gives the following generalization of Corollary~19:

\medskip

{\em Theorem~30: Second Order Generalized AEP:}
Let $\Xp$ be a stationary ergodic random field
with marginal distribution $P$ on $A$, and let 
$Q$ be an arbitrary probability measure on $\Ahat.$
Assume that the uniform $\alpha$-mixing coefficients
of $\Xp$ satisfy 
(\ref{eq:LIL-cond2})
and that $D_3=E_{P\times Q}[\rho^3(X,Y)]$ is
finite. Then for any $D\in(\Dmin,\Dav)$, and
with $g(x)$ defined as in (\ref{eq:functiong}),
$$
-\log Q^{n^d}(B(X_1^n,D))= n^dR_1(P,Q,D) + \sum_{u\in C(n)}g(X_u)
        + \frac{d}{2}\log n + O(\log\log n)
\;\;\;\;\mbox{w.p.1.}$$

\subsection{Applications}
Next we discuss applications of the second order 
generalized AEP to the 
$d$-dimensional analogs of the
data compression and pattern 
matching problems of Section~4.
As in Section~6.3, the only 
results stated explicitly are
those whose extensions to $\IN^d$
require modifications.

As mentioned in Section~6.3.1,
the one-dimensional construction 
of the random codes,
as well as the main tool used 
in their analysis, Theorem~9,
immediately generalize to the 
random field case. And since
all the second order results 
of Section~5.1 (Theorems~20--23)
are stated for $\iid$ sources, 
their statements as well as 
proofs carry over {\sl verbatim} 
to this case.

For the problem of waiting times,
we can use the second order generalized 
AEP of Theorem~30 to refine the SLLN
of Theorem~29
\ben
\frac{1}{n^d}\log W_n^d \to R_1(P,Q,D)
\;\;\;\;\mbox{w.p.1}
\een
to a corresponding CLT and LIL
as in the one-dimensional case.
These refinements are stated in
Theorem~31 below. Its proof is
identical to that of Theorem~24
in the one dimensional case. The
only difference here is that we 
need to invoke the CLT and LIL
for the partial sums of the random 
field $\{g(X_u)\;;\;u\in\IN^d\}$. 
Under the conditions of the 
theorem, these follow from the 
almost sure invariance principle
of \cite[Theorem~1]{berkes-morrow}.

\medskip

{\em Theorem~31:} 
Let $\Xp$ be a stationary ergodic random 
field and $\Yp$ be $\iid$, with marginal
distributions $P$ and $Q$ on $A$ and $\Ahat$,
respectively. Assume that the 
uniform $\alpha$-mixing
coefficients of $\Xp$ satisfy (\ref{eq:LIL-cond2})
and that $D_3=E_{P\times Q}[\rho^3(X,Y)]$ is
finite. Then for any $D\in(\Dmin,\Dav)$ the
following series is absolutely convergent
\be
\sigma^2\bydef 
\sum_{u\in\IN^d} E_P[g(X_{\orig})g(X_u)]
\label{eq:variance-d}
\ee
with $g(x)$ defined as in (\ref{eq:functiong}),
and, moreover:
\begin{itemize}
\item[]{\bf (CLT)} With $R_1=R_1(P,Q,D)$:
$$\frac{\log W^d_n \;-\; n^dR_1}{n^{d/2}}
        \weakly N(0,\sigma^2).$$
\item[]{\bf (LIL)}
The set of limit points of the sequence
$$\left\{
        \frac{\log W^d_n \;-\; n^dR_1}
             {\sqrt{2n^d\log\log n}}
  \right\},\quad n\geq 3$$
coincides with $[-\sigma,\sigma]$, with
probability one.
\end{itemize}

%
%
%

\section*{Acknowledgments}
We thank Tam\'{a}s Linder and Yuval Peres
for useful discussions regarding Theorems 7 and 8.

\appendix
\section{Proof of Theorem~7}
We prove the upper and lower bounds separately.
For the upper bound, 
recalling the definition
of $r_n(X_1^n)$ in
(\ref{eq:ratio})
we observe that
$$r_n(X_1^n,D)
\leq \frac{1}{n}\log {P_n(B(X_1^n,D))} -\frac{1}{n}\log Q^n(X_1^n)$$
where the second term converges to $H(P) + H(P\|Q)$ 
as $n\to \infty$, by the ergodic theorem. 
Since the first term is increasing in $D$,
for any fixed $D>0$ we have with 
$\BBP$-probability one:
\be
\limsupnd 
r_n(X_1^n,D)
	\leq 
	H(P) + H(P\|Q) + 
	\limsup_{n\to\infty}
	\frac{1}{n}\log P_n(B(X_1^n,D)).
\label{eq:discUB1}
\ee
Now the pointwise source coding
theorem (see \cite[Theorems~1 and~5]{konto-zhang:00})
implies that 
\be
\liminf_{n\to\infty}-\frac{1}{n}\log P_n(B(X_1^n,D))\geq R(D)
\;\;\;\;\mbox{w.p.1}
\label{eq:discUB2}
\ee
where $R(D)$ is the rate-distortion
function of the source $\Xp$
(in nats). 
%
From equations (\ref{eq:discUB1}) 
and (\ref{eq:discUB2}) we get
\ben
\limsupnd 
r_n(X_1^n,D)
& \leq & 
	H(P) + H(P\|Q) -R(D)\\
& \leq & 
	H(P) + H(P\|Q) - H(\BBP) + H(P) - R_1(D)
\;\;\;\;\mbox{w.p.1}
\een
where $R_1(D)$ denotes the first order 
rate-distortion function of $\Xp$,
$H(\BBP)$ is the entropy rate of $\Xp$
(both in nats), and
the second inequality follows
from the Wyner-Ziv bound;
see \cite[Remark~4]{wyner-ziv:71}.
The assumption that $\rho(x,y)=0$
if and only if $x=y$ implies that
$\lim_{D\to 0} R_1(D)=H(P)$,
so letting $D\downarrow 0$ 
the above right hand side becomes
$H(P) + H(P\|Q) -H(\BBP)$
and it is an easy calculation 
to verify that this is 
indeed the same
as $H(\BBP\|\BBQ)$.
This gives the required
upped bound.

For the lower bound we proceed
similarly by noting that 
$$
r_n(X_1^n,D)
\geq \frac{1}{n}\log {P_n(X_1^n)} 
	-\frac{1}{n}\log Q^n(B(X_1^n,D)),$$
where the first term converges to $H(\BBP)$
by the classical AEP
(as $n\to \infty$).
Since the second term is decreasing 
in $D$, for any fixed $D>0$ small
enough we have
with probability one:
\ben
\liminfnd 
r_n(X_1^n,D)
& \geq & 
        - H(\BBP) -
        \limsup_{n\to\infty}
        \frac{1}{n}\log Q^n(B(X_1^n,D))\\
& = & 
	- H(\BBP) + R_1(P,Q,D)
\een
where the last step follows from 
the generalized AEP in Theorem~1
(note that $\Dmin=0$ here).
By the characterization of the
rate-function in Proposition~2
we know that 
$$R_1(P,Q,D) = \sup_{\la'\leq 0} [\la' D-\LA(\la')]
	\geq [\la D-\LA(\la)]= 
	-E_{P}\left[\log E_{Q}\left(
		e^{\lambda(\rho(X,Y)-D)}
	\right)\right]$$
for any fixed $\la<0$.
Therefore, for any
$D$ small enough and
$\la<0$ we have
\ben
\liminfnd
r_n(X_1^n,D)
\geq - H(\BBP)  -E_{P}\left[\log E_{Q}\left(
                e^{\lambda(\rho(X,Y)-D)}
        \right)\right]
\;\;\;\;\mbox{w.p.1.}
\een
Letting $D\to 0$ and then $\la\to-\infty$,
by the dominated convergence theorem (and
the assumption $\rho(x,y)=0$ iff $x=y$)
the right hand side above converges
to $- H(\BBP) + H(P\|Q) + H(P)
= H(\BBP\|\BBQ),$ proving the
lower bound. 

Finally, since for each fixed $n$ 
the limit as $D\downarrow 0$ of
$r_n(X_1^n,D)$
exists,
it follows that 
the repeated limit
$\lim_{n}\lim_{D}$
also exists and is equal
to the double limit $H(P\|Q)$.
\qed

\section{Proof of Theorem~8}
Part~(a): Fixing $n$, let $f_n=dP_n/dQ_n$ and consider the set
$$
A_n \bydef \left\{ x_1^n : \; Q_n(B(x_1^n,D))>0 \;\; \forall D>0,
f_n(x_1^n) = \limsup_{D\downarrow 0}\; \frac{P_n(B(x_1^n,D))}{Q_n(B(x_1^n,D))}
= \liminf_{D\downarrow 0}\; \frac{P_n(B(x_1^n,D))}{Q_n(B(x_1^n,D))}
\, \right\}.
$$
By the Radon-Nikodym theorem 
(cf. \cite[Theorems 1.6.1, 1.6.2]{evans-gariepy}),
we know that $Q_n(A_n)=1$, hence also $P_n(A_n)=1$.
With $\BBP(\cup_n A_n^c)=0$, we conclude the proof of part~(a)
by applying Theorem~6 for $M_n=Q^n$ (in which case $H_n \geq 0$).

Part~(b): As $Q(A_1)=1$, in particular
$Q(B(x,D))>0$ for all $D>0$ and $Q$-almost
every $x\in\RL^d$ (hence also for $P=P_1$-almost
every $x\in\RL^d$), implying that $\Dmin$ of 
(\ref{eq:Dmin}) is zero. The same argument 
yields also that $P(B(x,D))>0$ for all $D>0$ 
and $P$-almost every $x$, hence $\Dmin$ is
still zero if we replace $Q$ by $P$. Thus, 
for all 
$D< \min\{E_{P\times Q}[\rho(X,Y)],E_{P\times P}[\rho(X,Y)]\}$,
applying Theorem~1 twice we get 
$$
  \lim_{n\to\infty}\; 
r_n(X_1^n,D)
	= R_1(P,Q,D)-R_1(P,P,D)
	\;\;\;\;\mbox{w.p.1.}
$$
For any probability measure $\mu$ and any $\la \leq 0$, let
$$
\Lambda(\lambda;\mu) = \int  \left[ 
\log \int e^{\lambda \rho(x,y)} d\mu(y)\right] dP(x).
$$

Fixing $D>0$ small enough,
we have by Proposition~2 
that $R_1(P,P,D) = \lambda D -\Lambda(\lambda;P)$
for the unique $\lambda=\lambda(D)<0$ such that $\Lambda'(\lambda;P)=D$, 
whereas $R_1(P,Q,D) \geq \lambda D - \Lambda(\lambda;Q)$. 
Since $E_{P\times P}[\rho(X,Y)]>0$, we have also that 
$\lambda(D) \downarrow -\infty$ as $D \downarrow 0$ (see (\ref{eq:la-lim})).
Consequently, 
$$
\liminf_{D\downarrow 0} \{ R_1(P,Q,D)-R_1(P,P,D) \} \geq 
\liminf_{\lambda \downarrow -\infty}
\{ \Lambda(\lambda;P) - \Lambda(\lambda;Q) \}
$$
Similarly, by Proposition~2 we have
$R_1(P,Q,D) = \widetilde{\lambda} D -\Lambda(\widetilde{\lambda};Q)$
for $\widetilde{\lambda}<0$ such that $\Lambda'(\widetilde{\lambda};Q)=D$,
$R_1(P,P,D) \geq \widetilde{\lambda}D -
\Lambda(\widetilde{\lambda};P)$, and with $E_{P\times Q}[\rho(X,Y)]>0$,
also $\widetilde{\lambda}\downarrow -\infty$ when $D \downarrow 0$.
Therefore, it suffices to show that
\be
\label{eq:dn-lim}
\lim_{\lambda \downarrow -\infty}
\{ \Lambda(\lambda;P) - \Lambda(\lambda;Q) \} = H(P\|Q) \;.
\ee

To this end, for any $\lambda <0$ and $x \in \RL^d$, let
$$
h_\lambda(x) \bydef
 \frac{E_P(e^{\lambda \rho(x,Y)})}{E_Q(e^{\lambda \rho(x,Y)})}
$$
noting that
$$
\Lambda(\lambda;P) - \Lambda(\lambda;Q) = \int \log h_\lambda(x) dP(x).
$$ 
Using the change of variable $U=\rho(x,Y) \geq 0$ followed
by integration by parts, we see that
$$
h_\lambda(x) =
 \frac{\int_0^\infty e^{\lambda u} g_{x}(u) du}{\int_0^\infty
e^{\lambda u} k_{x}(u) du} \;,
$$
where $g_x (r)=P(B(x,r))$ and $k_x(r)=Q(B(x,r))$ are nonnegative,
nondecreasing and bounded above by $1$. Considering separately
$u \leq 2\eta$ and $u>2\eta$, it is easy to check that for any $\eta>0$,
\be
\sup_{0 <r \leq 2\eta} \frac{g_{x}(r)}{k_{x}(r)} + \psi_{\lambda,x} \geq
h_\lambda(x) \geq \inf_{0 < r \leq 2\eta} \frac{g_{x}(r)}{k_{x}(r)}
\: \frac{1}{1+\psi_{\lambda,x}}
\label{eq:bd-dn}
\ee
where
\be
\psi_{\lambda,x} \bydef
\frac{\int_{2\eta}^\infty e^{\lambda u} du}
{\int_0^{2\eta} e^{\lambda u} k_{x}(u) du} \leq
\frac{1}{\eta |\lambda| k_x(\eta)} \;.
\label{eq:bd2-dn}
\ee
Fix $x \in A_1$ of part (a), in which case $k_x(r)>0$ for all $r>0$ and
$g_x(r)/k_x(r) \to f_1(x)$ as $r \to 0$.
Letting $\lambda \downarrow -\infty$ and then $\eta \to 0$, it
follows by (\ref{eq:bd-dn}) and (\ref{eq:bd2-dn}) that
$$
\lim_{\lambda \downarrow -\infty} h_\lambda(x) = f_1(x) \,.
$$
Recall that $P(A_1)=1$ and our assumption that
$\int \log k_x(\eta) dP(x) > -\infty$ for any $\eta >0$.
By our integrability conditions, the function
$\min\{0, \inf_{\lambda \geq 1} \log h_\lambda(x)\}$ is $P$-integrable,
hence, by Fatou's lemma,
$$
\liminf_{\lambda \downarrow -\infty} \int \log h_\lambda(x) dP(x) \geq 
\int \log f_1(x) dP(x) = H(P\|Q) \,.
$$
Moreover, in case $H(P\|Q)<\infty$, our assupmtions imply that
$\sup_{\lambda \geq 1} |\log h_\lambda(x)|$ is $P$-integrable,
hence by dominated convergence,
$\int \log h_\lambda(x) dP(x) \to \int \log f_1(x) dP(x)$ for
$\lambda \downarrow -\infty$,
as required to complete the proof of
(\ref{eq:dn-lim}).
\qed

\section{Proof of Theorem~27}

Recall our assumption that, for 
$\BBP$-a.e. $x_{[0,\infty)]}$, conditional on 
$X_{[0,\infty)]}=
x_{[0,\infty)]}$ 
the random variables $\{\rho_n(x_{C(n)},Y_{C(n)})\}$ satisfy the LDP
with a {\it deterministic} convex good rate-function 
denoted hereafter $R(\BBP,\BBQ,\cdot)$. Since 
$\rho$ is bounded, by Varadhan's lemma and convex duality, 
this implies that 
\be
\label{eq:am-las}
R(\BBP,\BBQ,D) =
\sup_{\lambda \in \RL} [ \lambda D - \Lambda_\infty(\lambda) ]
\bydef
\Lambda_\infty^*(D)
\ee
where for any $\lambda \in \RL$, the finite, deterministic limit
$$
\Lambda_\infty(\lambda) \bydef \lim_{n \to \infty} 
\frac{1}{n^d} \log \int e^{\lambda \sum_{u \in C(n)}\rho(x_u,y_u)} 
dQ_n(y_{C(n)})
$$
exists for $\BBP$-a.e. $x_{[0,\infty)}$ 
(cf. \cite[Theorem 4.5.10]{dembo-zeitouni:book}).
By bounded convergence, 
$\Lambda_\infty(\lambda)$ is also the limit of 
$$
\Lambda_n(\lambda) \bydef
\frac{1}{n^d} \int \left[
\log \int e^{\lambda \sum_{u \in C(n)}\rho(x_u,y_u)} 
dQ_n(y_{C(n)}) \right] dP_n(x_{C(n)}) \;.
$$

By stationarity, 
\be
\label{eq:am-dav}
\Dav=E_{P_n \times Q_n} (\rho_n(X_{C(n)},Y_{C(n)})), 
\;\;\;\forall n \geq 1
\ee
so replacing $P_1$, $Q_1$ and $\rho(x,y)$ of Proposition~2 by
$P_n$, $Q_n$ and $n^d \rho_n(x_{C(n)},y_{C(n)})$, 
respectively, we see that
\be
\label{eq:am-lasn}
R_n(P_n,Q_n,D) =
\sup_{\lambda \in \RL} [ \lambda D - \Lambda_n(\lambda) ]
\bydef
\Lambda_n^*(D) \,.
\ee
Note that 
$|\Lambda_n(\lambda)-\Lambda_n(\lambda')| \leq c |\lambda-\lambda'|$ for some
$c<\infty$ and all $n$, $\lambda,\lambda' \in \RL$,
hence the convergence of $\Lambda_n(\cdot)$ to $\Lambda_\infty(\cdot)$ is 
uniform on compact subsets of $\RL$. In particular, the convex,
continuous functions $\Lambda_n(\cdot)$ converge infimally to $\Lambda_\infty(\cdot)$,
and consequently, by \cite[Theorem 5]{wijsman}, the convex functions
$\Lambda_n^*(\cdot)$ converge infimally to $\Lambda_\infty^*(\cdot)$, that is
\be
\label{eq:inf-conv}
\Lambda_\infty^*(D)
=
\lim_{\delta \to 0} \limsup_{n \to \infty} \inf_{|\hat{D}-D|<\delta} \Lambda_n^*(\hat{D}) 
=
\lim_{\delta \to 0} \liminf_{n \to \infty} \inf_{|\hat{D}-D|<\delta} \Lambda_n^*(\hat{D}) \,.
\ee

It follows from (\ref{eq:am-dav}) and Jensen's inequality 
that $\Lambda_n(\lambda) \geq \lambda \Dav$ for all $n$ and $\lambda$,
hence, for $D \leq \Dav$ suffices to consider $\lambda \leq 0$ in
(\ref{eq:am-las}) and in (\ref{eq:am-lasn}). Thus, for $1 \leq n \leq \infty$,
$\Lambda^*_n$ are non-negative, convex,
and monotone non-increasing on $[0,\Dav]$, with 
$\Lambda^*_n(\Dav)=0$. For $1 \leq n \leq \infty$,
let
$$
\Dminn \bydef \lim_{\lambda \downarrow -\infty} \frac{\Lambda_n(\lambda)}{\lambda} \,,
$$
so that $\Lambda_n^*(D)=\infty$ for $D < \Dminn$, while
$\Lambda_n^*(D)<\infty$ for $D>\Dminn$.
Note that for $n < \infty$ this coincides with the definition of $\Dminn$ 
given in (\ref{eq:dminn}). It is easy to check then that (\ref{eq:inf-conv})
implies the pointwise convergence of $\Lambda^*_n(\cdot)=R_n(\BBP,\BBQ,\cdot)$
to $\Lambda^*_\infty(\cdot)=R(\BBP,\BBQ,\cdot)$ at any $D$ for which
$\Lambda^*_\infty(D-\delta) \downarrow \Lambda^*_\infty (D)$, that is,
for all $D \neq \Dinf$. In particular, necessarily $\Dinf \in [\Dmin,\Dav]$,
and $\Dinf$ may also be defined via (\ref{eq:dinf}). 
The continuity of $R(\BBP,\BBQ,D)$ at $D \in (\Dmin,\Dav)$, $D \neq \Dinf$ implies
the equality in (\ref{eq:ldp-27}) for such $D$, thus completing the proof of the 
theorem.
\qed

\section{Proof of Theorem~28}
For each $m\geq 1$, let $G_m$ be the 
collection of ``good'' realizations
$x_{\IN^d}\in A^{\IN^d}$
$$G_m= \left\{ x_{\IN^d}\in A^{\IN^d} 
: \;\;Q_n(B(x_{C(n)},D))>0 \;\mbox{for all}\; n\geq m \right\}$$
so that the assumption that 
$Q_n(B(X_{C(n)},D))>0$ eventually, with probability one 
translates to 
\be
\BBP\{\cup_{m\geq 1} G_m\}=1.
\label{eq:eventuallyOK}
\ee

To prove the lower bound
we choose and fix an $m\geq 1$
and a realization $x_{\IN^d}\in G_m$. 
Then for any $K>1$:
\ben
\PR\{W_n^d<K\,|\,X_{C(n)}=x_{C(n)}\}
&\leq& \sum_{u\in[0,\lfloor K^{1/d} \rfloor-1]^d}
	\,Q_n\{Y_{u+C(n)}\in B(x_{C(n)},D)\}\\
&\leq& K\,Q_n(B(x_{C(n)},D)).
\een
Since, by its definition,
$W_n$ is always greater than
or equal to one, this 
inequality trivially holds 
also for $K\in(0,1]$.
Setting 
$K=[n^{1+\epsilon}Q_n(B(x_{C(n)},D))]^{-1}$ 
above gives,
for all $n \geq m$,
\ben
\PR\{\log[W^d_nQ_n(B(X_{C(n)},D))]<-(1+\epsilon)\log n
\,|\,X_{C(n)}=x_{C(n)}\}
\leq \frac{1}{n^{1+\epsilon}}.
\een
Since this bound is uniform over 
$x_{\IN^d} \in G_m$ and summable, the Borel-Cantelli 
lemma and assumption (\ref{eq:eventuallyOK})
imply that
\be
\log[W^d_nQ_n(B(X_{C(n)},D))]\;\geq\;-(1+\epsilon)\log n
\;\;\;\;\mbox{eventually, w.p.1.}
\label{eq:bc:2}
\ee

For the upper bound 
we choose and fix an $m\geq 1$
and a realization $x_{\IN^d}\in G_m$,
and take $K\geq (n+1)^d$.
Note that 
$$\PR\{W_n^d>K\,|\,X_{C(n)}=x_{C(n)}\}
\leq\Pr
    \left\{
      \sum_{u\in [0,M]^d}
      \IND_{
	\{
	Y_{nu+C(n)}\in B(X_{C(n)},D)
	\}
	   }
		= 0 
    \right\}
$$
where 
the sum is over the $(M+1)^d$
integer positions $u\in[0,M]^d\subset\IN^d$,
$nu$ denotes the point 
$(nu_1,nu_2,\ldots,nu_d)\in\IN^d$,
and
$$M=M(K,n)\bydef
\left\lfloor\frac{K^{1/d}-1}{n}
\right\rfloor.$$
Let $\Sigma_n$ denote the sum 
in the above probability,
$$\Sigma_n=\sum_{u\in[0,M]^d}I_n(u)$$
where $I_n(u)$ is the indicator function 
of the event $\{Y_{nu+C(n)}\in B(X_{C(n)},D)\}$.
In this notation:
\be
\PR\{W_n^d>K\,|\,X_{C(n)}=x_{C(n)}\}
\,\leq\,\BBQ\{\Sigma_n=0\}\,\leq\,\frac{\VAR_\BBQ(\Sigma_n)}
				{[E_\BBQ(\Sigma_n)]^2}.
\label{eq:estimate1}
\ee
By stationarity
\be
E_\BBQ(\Sigma_n) = [M+1]^dQ_n(B(x_{C(n)},D))
\label{eq:estimate2}
\ee
and by the definition of the
$\phi$-mixing coefficients, if $u\neq v$,
$$E_\BBQ\{I_n(u)I_n(v)\}\leq Q_n(B(x_{C(n)},D))
[\phi_n(nd(u,v)-n+1)+Q_n(B(x_{C(n)},D))].$$ 
Using the last two estimates
we can bound the variance as
\be
\VAR_\BBQ\{\Sigma_n\}
&=&\sum_{u,v\in[0,M]^d}
\COV_\BBQ(I_n(u), I_n(v))
        \nonumber\\
&\leq&[M+1]^dQ_n(B(x_{C(n)},D))
	\nonumber\\
& & \quad
	+\sum_{u,v\in[0,M]^d,\;u\neq v}
	\Big[Q_n(B(x_{C(n)},D))
	\phi_n(nd(u,v)-n+1)\Big]
	\nonumber\\
&\leq&[M+1]^dQ_n(B(x_{C(n)},D))
	\left[1+\sum_{j=1}^{M}
c_d j^{d-1}\phi_n(nj-n+1)\right]
        \label{eq:estimate3}
\ee
where
$c_d j^{d-1}$ bounds the number of possible points
$u$ that can be at a distance exactly $j$ from a 
given point $v$ (for some constant $c_d$).
By assumption (\ref{eq:dobrushin}) we can find
a finite constant $\Phi$ such that the expression
in square brackets in (\ref{eq:estimate3}) is bounded
above by $\Phi$, uniformly in $n$.
Substituting this bound, together with
(\ref{eq:estimate2}) and (\ref{eq:estimate3}),
in (\ref{eq:estimate1}), gives
\be
\PR\{W_n>K\,|\,X_{C(n)}=x_{C(n)}\}&\leq&\frac
        {\Phi}
        {[M+1]^dQ_n(B(x_{C(n)},D))}.
\label{eq:estimate4}
\ee
Let $\epsilon>0$ arbitrary, 
take $n$ large enough so 
that $n^{(1+\epsilon)/d}\geq 2$,
and let $K=n^{d+1+\epsilon}/Q_n(B(x_{C(n)},D)).$
Simple algebra shows that with this choice
of $K$ we have
$$[M+1]^dQ_n(B(x_{C(n)},D))\geq\frac{1}{2}n^{1+\epsilon}$$
and substituting this in (\ref{eq:estimate4}) 
yields
\ben
	\PR\{\log[W_n^dQ_n(B(
x_{
C(n)},D))] >
	(d+1+\epsilon)\log n\,|\,X_{C(n)}=x_{C(n)}\}\leq
        \frac{2\Phi}{n^{1+\epsilon}}.
\een
This bound is uniform 
over $x_{\IN^d} \in G_m$ and summable, 
so the Borel-Cantelli lemma 
and (\ref{eq:eventuallyOK})
imply that 
\be
\log[W_n^dQ_n(B(
X_{
C(n)},D))]\;\leq\;(d+1+\epsilon)\log n
\;\;\;\;\mbox{eventually, w.p.1.}
\label{eq:bc:1}
\ee

Combining (\ref{eq:bc:1}) and (\ref{eq:bc:2}) 
completes the proof.
\qed


\end{document}